\documentclass[graybox]{svmult}

\usepackage{type1cm}        
\UseRawInputEncoding                          
\usepackage{makeidx}         
\usepackage{graphicx}        
\usepackage{multicol}        
\usepackage[bottom]{footmisc}

\usepackage{newtxtext}       %
\usepackage{newtxmath}       

\makeindex             

\definecolor{Green}{rgb}{0,.8,0}
\definecolor{Green}{rgb}{0.20,0.43,0.09}  

\newcommand{\C}{\mathbb{C}}

\def\CC{{\cal C}}

\newcommand{\hlie}{\mathfrak{h}}

\newcommand{\adiag}{{A_{\mathrm{diag}}}}
\def\hsub{\mathfrak{h}_{sub}}
\usepackage[colorlinks=true, linkcolor=blue, 
hyperfootnotes=true,citecolor=blue,urlcolor=black]{hyperref}

\begin{document}

\title*{Differential Galois Theory and Integration}
\author{Thomas Dreyfus and Jacques-Arthur Weil}
\institute{Thomas Dreyfus \at Institut de Recherche Math\'ematique Avanc\'ee, U.M.R. 7501 Universit\'e de Strasbourg et C.N.R.S. 7, rue Ren\'e Descartes 67084 Strasbourg, France. \email{ dreyfus@math.unistra.fr}
\and Jacques-Arthur Weil \at XLIM,  U.M.R 7252, Universit\'e de Limoges et C.N.R.S, 123 avenue Albert Thomas, 87060 Limoges Cedex, France. \email{Jacques-Arthur.Weil@unilim.fr}}

%
%
\maketitle

\abstract{In this chapter, we present methods to simplify reducible linear differential systems before solving. Classical integrals appear naturally as solutions of such systems. We will  illustrate the methods developed in \cite{DrWe19a} on several examples to reduce the differential system.  This will give information on 
	potential algebraic relations between integrals.}
\keywords{
Ordinary Differential Equations, Differential Galois Theory, Computer Algebra, Integrals, Lie Algebras, $D$-finite functions.
}

\section*{Introduction}


In this chapter, we will review properties of block triangular linear differential systems and their use to compute properties of integrals. 
\par  
Let $\mathbf{k}=\C(x)$ and  $A \in \mathrm{Mat}(n,\mathbf{k})$. We will study the corresponding \emph{linear differential system}
$[A] : \, \partial_{x} Y=A Y$.
More generally, we might consider linear differential systems over a differential field $(\mathbf{k}, \partial)$ of characteristic zero, that is a field $\mathbf{k}$ equipped with an additive morphism $\partial$ satisfying the Leibniz rule $\partial (ab)=a\partial (b)+\partial(a)b$.\par 
The Galois theory of linear differential equations aims at understanding what are the algebraic relations between the solutions of $[A]$. We attach to $[A]$ a group that measures the relations. The computation of 
this differential Galois group is a hard task in full generality.  
The goal of this chapter is to illustrate on  examples the method described in \cite{DrWe19a} that focuses on the reduction of block-triangular linear differential systems. This approach is powerful enough to understand the desired relations  on the solutions. \par 

Given an invertible matrix $P\in \mathrm{GL}(n, \mathbf{k})$, the linear change of variables 
$Y=PZ$ produces a new differential system denoted $Z'=P[A]Z$
where
\[ P[A] := P^{-1} A P - P^{-1} P'.\]
Two linear differential systems $[A]$ and $[B]$ are called \emph{(gauge) equivalent} over $\mathbf{k}$ if there exists a gauge transformation, an invertible matrix $P\in \mathrm{GL}(n, \mathbf{k})$, such that $B=P[A]$. 
A linear differential system is called $\emph{reducible}$ (over $\mathbf{k}$) when it is gauge equivalent to a linear differential system in block triangular form:
\begin{equation}\label{redsys}
[A] : \, \partial Y=A Y, 
\, \textrm{ with } 
A=\left(\begin{array}{c|c}  {A_{1}} & 0 \\\hline  {S} &  {A_{2}}\end{array}\right) \in \mathrm{Mat}(\mathbf{k}).
\end{equation}

It turns out that computing properties of integrals of $D$-finite functions\footnote{A function is \emph{$D$-finite} when its is a solution of a linear differential equations with coefficients in $\mathbf{k}$} or of some types of iterated integrals may be reduced to computing solutions of such reducible linear differential systems. 
Computing the differential Galois groups of block triangular systems gives, in turn, information on properties of their solutions. This idea was promoted by Bertrand \cite{Be01a}  or Berman and Singer \cite{BeSi99a} who showed how to compute Galois groups of some reducible systems and how this would reveal algebraic properties of integrals.

Our aim, in this chapter, is to show how to compute such algebraic properties. The underlying theory is developed in \cite{DrWe19a} and \cite{CaWe18a};  general references for differential Galois theory are for example \cite{PuSi03a, CrHa11a, Si09a}; general references for constructive theory of reduced forms of differential systems are \cite{ApCoWe13a,MiSi02a, ApDrWe16a,BaClDiWe20a,BaClDiWe16a}.
We will rely on many examples rather than on cumbersome theory and provide references for interested readers.
For example, we will consider the dilogarithm function $\textrm{Li}_2$ defined by 
$\textrm{Li}_2'(x) = -\frac{\ln(1-x)}{x}$ and our $4$-dimensional example (in the next sections) will provide a simple algorithmic proof that it is not only transcendent but algebraically independent of $e^x, \ln(x)$ and $\ln(1-x)$; the proof will only use rational solutions of linear first order differential equations. \par 
The chapter is organized as follows. We begin by some examples to illustrate what  the method will provide. Then we {give} a review of differential Galois theory and reduced forms of differential systems. We finish by explaining the strategy of \cite{DrWe19a} on two examples that are chosen so that almost all calculations can be reproduced easily.

\textbf{Acknowledgment} This project has received funding from the ANR project \emph{De Rerum Natura} ANR-19-CE40-0018. We warmly thank the organizers of the conference \emph{Antidifferentiation and the Calculation of Feynman Amplitudes} for  stimulating exchanges and lectures.  We specially thank the referee for many observations and precisions that have enhanced the quality of this text. 

\section{Examples}\label{sec:1}
\subsection{A first toy example}
We consider two confluent Heun\footnote{See \url{https://dlmf.nist.gov/31.12} or \cite{Mo18p}. The notation $HeunC$ is the syntax in \textsc{Maple}.} functions
\begin{eqnarray*} 
	f_1(x) &=& \exp \left(\frac{\sqrt{3}}{12} x\right)
	\textrm{HeunC} \left( \frac{\sqrt {3}}{6},-\frac{1}{3},-\frac{1}{3},\frac{1}{48},\frac{11}{48}; x \right) 
		\\ &= & 
		1-{\frac{7}{96}}x-{\frac{719}{46080}}{x}^{2}-{\frac{127307}{10616832}
}{x}^{3}-{\frac{82319293}{10192158720}}{x}^{4}+O \left( {x}^{5} \right)
\end{eqnarray*} 
and 
\begin{eqnarray*}  
	f_2(x) &=& \exp\left(\frac{\sqrt{3}}{12} x\right) \sqrt [3]{x} \textrm{HeunC} \left(\frac{\sqrt {3}}{6}\,,\frac{1}{3},-\frac{1}{3},\frac{1}{48},{\frac{11}{48}} ; x \right)
	\\ & =& 
	x^{\frac13}
	\left(1+{\frac{25}{192}}x+{\frac{8977}{129024}}{x}^{2}+{\frac{1099183}{
26542080}}{x}^{3}+O ( {x}^{4} ) \right).
\end{eqnarray*}
They form a basis of solutions of the second order linear differential equation 
\[L(y):={\frac {
{\rm d}^{2}}{{\rm d}{x}^{2}}}y \left( x \right) 
 + \frac23\left( \frac1x+  \frac{1}{x-1} \right) {\frac {\rm d}{{\rm d}x}}y \left( x \right)
Ê-{\frac { \left( 3\,{x}^{2}-6\,x+7 \right) }{144\,x
 \left( x-1 \right) }}y \left( x \right) =0 
\]
The Wronskian relation gives us the algebraic relation $(f_1 f_2'- f_1' f_2)^3= x^2(x-1)^2$. 
The equation has order two so the Kovacic algorithm
\cite{Kov86,UlWe96a,HoWe05a} can be used to compute 
the differential Galois group
and we find that no other algebraic relations exist between 
$f_1$, $f_2$, $f_1'$ and $f_2'$. 
Now let $F_i(x):=\int^x f_i(t)dt$ be a primitive. We want to determine whether $F_1$ and $F_2$ are algebraically independent of the $f_i$ and $f_i'$ or not.
The techniques explained below will show that this question reduces to asking whether there is a \emph{rational} solution to the linear differential system
\[ Z' =  \left( \begin {array}{cc} 0&{\frac {-3\,{x}^{2}+6\,x-7}{144\,x
 \left( x-1 \right) }}\\ \noalign{\medskip}-1&\frac13\,{\frac {4\,x-2}{x
 \left( x-1 \right) }}\end {array} \right) \cdot Z 
 \, + \, \begin{pmatrix} 1 \\ 0 \end{pmatrix}
 \]
or, equivalently, whether the following linear differential equation (the left hand side turns out to  be the adjoint operator of $L$) has a rational solution:
\[ -g''(x) +\frac23\,\frac{ \left( 2\,x-1 \right)} {x \left( x-1 \right) } g'(x) +{\frac { \left( 3\,{x}^{4}-9\,{x}^{3}-179\,{
x}^{2}+185\,x-96 \right)  }{144\,{x}^{2} \left( x-1
 \right) ^{2}}} g \left( x \right)=1.
\]
It can be seen directly or with a computer algebra system that this equation has no  rational solution. The underlying theoretical tools are from the constructive differential Galois theory. 
However, in operational  terms, it is rather easy to compute and check : no hard theory is required for the calculation. Let us unveil a corner of the underlying tools. 

We have in fact studied the differential Galois group of the differential system $[A]$ with
\[ 
A=  \left( \begin {array}{ccc} 0&1&0\\  {\frac {3\,{x}^{
2}-6\,x+7}{144\,  (x-1)^3 x^2}}
&-\frac{2}{3}\left(\frac{1}{x-1} +  \frac{1}{x} \right)
&0\\  1&0&0\end {array}
 \right)   
 \]
which admits the fundamental solution matrix
 \[
U:= \left(\begin{matrix} 
       f_1(x) & f_2(x) & 0 \\
        f_1'(x) & f_2'(x) & 0 \\
        \int^x f_1(t) {\rm dt} & \int^x f_2(t) {\rm dt} & 1
    \end{matrix} \right).
\]
Our calculation of rational solution above shows (with the tools displayed below) that the differential Galois group of the system $[A]$ will have dimension $5$. This in turn shows that the integrals 
$ \int^x f_j(t) {\rm dt}$ are algebraically independent of $f_1$, $f_2$ and their derivatives. 
In fact, it even shows that both integrals are algebraically independent.

\subsection{A second toy example}
Recall that the hypergeometric function is given by the 
formula 
$${}_2F{}_1([a, b], [c])(x)=\displaystyle \sum_{n=0}^{\infty}\frac{(a)_{n}(b)_{n}}{(c)_{n}}\frac{x^{n}}{n!}, \quad \hbox{with } \quad (a)_0=1, \; (a)_n= a(a+1)\dots (a+n-1). $$
We now take two hypergeometric functions (with nice modular properties)
	$${f_1(x) = {}_2F_1\left(\left[-\frac{1}{3}, \frac{1}{12}\right], \left[\frac{7}{12}\right]\right)(x)}$$ and 
	$$f_2(x) = x^{5/12}{ }_2F_1 \left(\left[ \frac{1}{12},  \frac{1}{2}\right], \left[\frac{17}{12}\right]\right)(x).$$
The vectors $Y_i:=(f_i, f_i')^T$ are solutions of $Y'=A_1 Y$ below. To study properties of their integrals, 
we set $Y:=(f, f', \int^x f(t) dt)^T$ and we have $Y'=AY$ where
\[ A_1 = \left( \begin {array}{cc} 0&1\\  \frac{1}{36}\,{\frac {1}{x
 \left( x-1 \right) }}&-{\frac {7}{12\,x}}-\frac{1}{6\left( x-1 \right)}\,  \end {array} \right)
\textrm{ and }
 A=\left(\begin{matrix} A_1 & 0 \\ 1 \,0 & 0\end{matrix}\right).
\]

{We will see in the sequel how we may find a suitable change of variables } 
\[ Q:= \left( \begin {array}{ccc} 1&0&0\\ \noalign{\medskip}0&1&0
\\ \noalign{\medskip}-{\frac{15}{44}}+{\frac {45\,x}{44}}&-{\frac {9\,
x \left( x-1 \right) }{11}}&1\end {array} \right)
\]
such that 
\[ 
Q[A] = \left(\begin{matrix} A_1 & 0 \\ 0\,0 & 0\end{matrix}\right).
\]
This shows that \[
\int^x \!f_i \left( t \right) \,{\rm dt}=
	- \frac {9}{11}\,x \left( x-1 \right) f_i'\left( x \right) 
	+ \frac{15}{44} \left( 3 \, x -1 \right)  f_i \left( x \right) 
	{+}{c_i }, \quad c_i\in \C,
\]
and the differential Galois group of $[A]$ has dimension 3. 
We note that none of these two examples is new and that efficient methods to handle these questions have been developed by Abramov and van Hoeij in \cite{AbHo99b,AbHo97t}.

\subsection{Integrals via reducible systems}
These two examples have shown how, by augmenting the dimension of a linear differential system, we can study integrals of its solutions. Conversely, block triangular systems give rise to integrals via variation of constants; we review this for clarification. 
For a factorized reducible system
 $[A]$ of the form \[A = \left(\begin{array}{c|c} A_{1} & 0 \\\hline S & A_{2}\end{array}\right)
 = A_{diag} + A_{sub}, \]
 with  $A_{diag}:=\left(\begin{array}{c|c} A_{1} & 0 \\\hline 0 & A_{2}\end{array}\right)$ and $A_{sub}:=\left(\begin{array}{c|c}0 & 0 \\\hline S & 0\end{array}\right)$,
we have a fundamental solution matrix of the form
\[ U =\left(\begin{array}{c|c} U_{1} & 0 \\\hline U_2 V & U_{2}\end{array}\right) 
	= \left(\begin{array}{c|c} U_{1} & 0 \\\hline 0 & U_{2}\end{array}\right) \left(\begin{array}{c|c} \mathrm{Id}_{\mathrm{n_{1}}} & 0 \\\hline  V & \mathrm{Id}_{\mathrm{n_{2}}}\end{array}\right). 
\]
Once $U_1$ and $U_2$ are known, $V$ is given by integrals : $\partial V=U_2^{-1} S U_1$. 
So it may seem that no further theory is required. However, $U_1$ and $U_2$ are fundamental solution matrices for systems $\partial U_i =A_i U_i$ so the relation $\partial V=U_2^{-1} S U_1$ may involve integrals of complicated $D$-finite functions. 
\\
Our approach will be to first ``reduce'' $S$ as much as possible, using manipulations with rational functions, to prepare the system for easier solving. In return, we will obtain algebraic information on all (possibly iterated) integrals that may occur in the  relation $\partial V=U_2^{-1} S U_1$.

We note, for the record, that factorized differential operators also correspond to block triangular differential systems.

\subsection{Example of situations involving reducible linear differential systems}

\subsubsection{Operators from statistical physics and combinatorics} 
In many models attached to statistical mechanics, quantities are expressed as multiple integrals depending on a parameter. They are generally holonomic in this parameter, meaning that they are $D$-finite, i.e they are solutions to linear differential operators (for example, the so-called Ising operators). A description of this setting may be found in the books of Baxter \cite{Ba82a} and McCoy \cite{Mc10a} (notably Chapters 10 and 12) or in the  {surveys} \cite{McMa12a,McAsBoHaMaOr11a} (and references therein).
Similarly, in combinatorics, sequences that satisfy recurrence relations may be studied through their generating series, which are often $D$-finite.  

Experimentally (see, for example, \cite{BoHaMaMcWeZe07a, BoBoHaMaWeZe09a, BoBoHaMaWeZe10a,BoBoHaHoMaWe11a, BoHaMaWe14a,BoHaMaWe15a} or \cite{HaKoMaZe16l,Ko13a,AbBoKoMa20q,AbBoKoMa18z}), it turns out that many differential operators coming from these processes are factored into a product of smaller order factors - and the corresponding companion systems are reducible. 

Other cases of reducible systems are the ones admitting reducible monodromies 
such as \cite{SaTa16e} or in the works of Kalmykov on Feynman integrals via Mellin-Barnes integrals \cite{Ka06g,ByKaKn10v, KaKn12o,KaKn17c,KaKn11a}; about differential equations for Feynman integrals, one may consult the works of Smirnov \cite{Sm12e,LeSmSm18g}.

Finally, we mention the paper \cite{AbBlRaSc14b} with other examples of integrals where the techniques presented below may offer alternative approaches for some of the computations. 
	
\subsubsection{Variational equations of nonlinear differential systems}
Another natural source of reducible systems is the old method of \emph{variational equations} of nonlinear differential systems along a given particular solution $\phi$.
One can form the linear differential equation describing perturbations along this solution $\phi$. The general principle is that obstructions to integrability of the nonlinear system can be read on this linear differential system. 
Ziglin \cite{Zi82a}  linked non-integrability to non-commutations in the monodromy group with concrete versions given e.g.  in \cite{ChRo91a} and \cite{Sa14a}. It was generalized in  the theory of Morales-Ruiz and Ramis and extended with Sim\'o \cite{MoRa01b,MRS} for Hamiltonian systems; they prove that a Hamiltonian system is completely integrable only if all its variational equations have a virtually abelian differential Galois group. 
Extensions to nonhamiltonian differential equations can be found in \cite{AyZu10a,CaWe18a}. 
Applications to specific problems have been occasions to establish effective criteria. For example:
the three body problem \cite{Ts01c,BoWe03a}, $n$-body problems \cite{MoSi09a,Co12b}, 
Hill systems (movements of the moon) \cite{MoSiSi05a,ApWe12b}
or a 
swinging Atwood machine \cite{PPRSSW10a}.

These variational equations can be written in the form of reducible linear differential systems of big dimension. The simplification techniques outlined below are hence particularly relevant to make computation on such systems practical (see \cite{ApDrWe16a}). 

We note that the Morales-Ramis theory has had a spectacular recent development, initiated in \cite{Mo20s} where this variational approach is applied to path integrals
thus establishing a beautiful and unexpected bridge with the previous subsection.  

Reducible linear differential systems appear naturally  in another type of perturbative approach: in $\epsilon$-expansions of solutions for perturbed systems $\frac{d}{dx}Y=B(x,\epsilon) Y$ like the ones that appear in works of J. Bl\"umlein, C. Raab, C. Schneider, J. Henn and others for example. 
\section{Reduced forms of linear differential systems}
\subsection{Ingredient \#1 : Differential Galois-Lie Algebra}
In what follows, $(\mathbf{k},\partial)$ is a differential field of characteristic $0$. We outline a brief exposition of the Galois theory of linear differential equations, see  \cite{PuSi03a, CrHa11a, Si09a} for expositions with proofs. 
We consider a linear differential system
\begin{equation}
[A]: \quad {\partial Y}=A Y, \quad A \in \mathrm{Mat}(n,\mathbf{k}).
 \end{equation}
Let $\CC$ be the field of constants of the differential field $\mathbf{k}$, that is $\CC:=\{a\in \mathbf{k}| \partial (a)=0\}$. We will assume that $\CC$ is algebraically closed, {i.e, every non constant polynomial equation has a solution in $\CC$}. 

A \emph{Picard-Vessiot  extension} is a field $K=\mathbf{k}(U)$, where $U$ is a fundamental solution matrix\footnote{An invertible $n\times n$ matrix $U$ such that 
$U'=AU$.} of $[A]$, such that the field of constants of $K$ is still $\CC$. 
This can be constructed algebraically; alternatively, when $\mathbf{k}$ is a field of meromorphic functions, one may consider a local matrix of power series solutions at a regular point and this gives a Picard-Vessiot extension. 
A Picard-Vessiot extension is unique modulo differential field isomorphisms.
\par 

The \emph{differential Galois group} $G:=\mathrm{Aut}_{\partial}(K/\mathbf{k})$ is the set of automorphisms of $K$ which leave the base field $\mathbf{k}$ fixed and commute with the derivation. 
Let 
$\sigma\in G$. By construction, $\sigma (U)$ is also a fundamental solution matrix in $K$ and
we find
that there exists a matrix $[\sigma]\in \mathrm{GL}(n, \CC)$ such that $\sigma (U)=U\cdot[\sigma]$. The map $\sigma \mapsto [\sigma]$ provides a faithful representation of $G$ as a subgroup of $\mathrm{GL}(n, \CC)$, actually a linear algebraic group. If we change the fundamental solution, we obtain a conjugate representation.\par 

We recall that, given an invertible matrix $P\in \mathrm{GL}(n, \mathbf{k})$, the linear change of variables $Y=PZ$ produces a new differential system denoted $\partial Z=P[A]Z$ where $$P[A] := P^{-1} A P - P^{-1} \partial P.$$
 We note that such a gauge transformation $P[A]$, with $P\in \mathrm{GL}(n, \mathbf{k})$, does not change the Galois group. \par 
Given a Picard-Vessiot  extension $K=\mathbf{k}(U)$, the polynomial relations  among all entries of $U$ (over $\mathbf{k}$) form an ideal $I$. 
The Galois group $G$ can then be viewed as the set of matrices which stabilize this ideal $I$ of relations. 
Thus the computation of $G$ is strongly related to the understanding of the algebraic relations among the solutions.

The \emph{Galois-Lie algebra of $[A]$} is the Lie algebra $\mathfrak{g}$ of the differential Galois group $G$. It is defined as the tangent space of $G$ at the identity $\mathrm{Id}$. The dimension of $\mathfrak{g}$ measures the transcendence degree of $K$ over $\mathbf{k}$, that is  
\[ \dim_{\CC}\mathfrak{g} ={\rm  trdeg}(K/\mathbf{k}).\]
One way of computing the Lie algebra (\cite{PuSi03a}) is the following : $\mathfrak{g}$ is the set of matrices $N$ such that $ \mathrm{Id} +\epsilon N \in G(\CC[\epsilon])$ with $\epsilon^2=0$. In other terms, $ \mathrm{Id} +\epsilon N$ satisfies the defining equations of the group modulo $\epsilon^2$. 
For example, 
${\mathrm{SL}(n,\CC)}$ (set of $M$ such that $\det(M)=1$) gives the Lie algebra $\mathfrak{sl}(n,\CC)$ of 
matrices $N$ such that $\mathrm{Tr}(N)=0$. The symplectic group ${\mathrm{Sp}(2n,\CC)}$ is the set of $M$ such that $M^{T}\cdot J \cdot M = J$; its Lie algebra $\mathfrak{sp}(n,\CC)$ is found to the set of $N$ such that $N^{T}J+JN=0$, with 
$J=\left(\begin{array}{cc}0& \mathrm{Id}  \\ -\mathrm{Id} & 0\end{array}\right)$. 
The additive group $\mathbb{G}_{a}:=\left\{\left(\begin{array}{cc}1 & a \\0 & 1\end{array}\right), a\in \CC\right\}$ admits the Lie algebra ${\mathfrak{g}_{a}}=Span_{\CC}\left\{ \left( \begin{array}{cc}0 & 1 \\0 & 0\end{array}\right) \right\}$; the multiplicative group $\mathbb{G}_{m}=\left\{\left(\begin{array}{cc}a & 0 \\0 & \frac1a\end{array}\right), a\in \CC^{*}\right\}$ admits the Lie algebra $\mathfrak{g}_{m} = Span_{\CC}\left\{ \left(\begin{array}{cc} 1 & 0 \\0 & -1\end{array}\right) \right\}$.

The reduction technique exposed in this chapter aims at computing directly the Galois-Lie algebra $\mathfrak{g}$ before computing $G$ itself. Although the theory is not obvious, the resulting calculations are reasonably simple. 

\subsection{Ingredient \#2 : Lie algebra $\textrm{Lie}(A)$ associated to $A$ }
The Lie algebra $\textrm{Lie}(A)$ associated to the matrix $A$ is defined as follows. 
Let  $a_1,\ldots, a_s\in \mathbf{k}$ be a basis of the $\CC$-vector space generated by the coefficients of $A$.  
We can then decompose $A$ as $A=\sum_{i=1}^s a_i M_i$ where the $M_i$ are constant matrices.
Now, we consider the smallest Lie algebra containing all the $M_i$: this is the vector space generated by the $M_i$ and all their iterated Lie brackets ($[M,N]:=MN-NM$); then we take its algebraic envelope. 
\begin{definition}
The  Lie algebra $\textrm{Lie}(A)$ associated to the matrix $A$ is the smallest \emph{algebraic}\footnote{A Lie algebra is called algebraic when it is the Lie algebra of a linear algebraic group. When the $M_i$ are given, this can be computed, see \cite{FiGr07a}, Section 3, or  \cite{DrWe19a}, Section 6.}
 Lie algebra containing  
 all the $M_i$.
 \end{definition}
{The decomposition  $A=\sum_{i=1}^s a_i M_i$ is not unique but the vector space generated by the $M_i$ is unique. Thus, the associated Lie algebra $\textrm{Lie}(A)$ does not depend on the chosen decomposition. }
 \\
 This Lie algebra $\textrm{Lie}(A)$ appears in works of Magnus \cite{Ma54a} or Feynman who use the Baker-Campbell-Hausdorf formula to write solutions of $\partial Y=AY$ as (infinite) products of exponentials constructed with Lie brackets. Wei and Norman give in \cite{WeNo63a,WeNo64a} a finite formula to solve the system when $\textrm{Lie}(A)$ is solvable. This formula is well-known in physics and control theory but not as well among mathematicians. 
 The terminology of Lie algebra associated to $A$ appears\footnote{For this  reason, some authors, including ourselves, call the decomposition $A=\sum_{i=1}^s a_i M_i$ a \emph{Wei-Norman decomposition of $A$}.} in \cite{WeNo63a,WeNo64a} (in there, it is defined as the Lie algebra generated by all values of $A(z_0)$ for $z_0$ spanning all constants minus singularities and the algebraic envelope is missing). In the sequel, we will study a $4$-dimensional example and an $8$-dimensional example where our technique will have some relations to the Wei-Norman approach; namely, we will change $A$ to obtain an associated Lie algebra $\textrm{Lie}(A)$ of minimal dimension so that solving formulas become optimal in some sense.
 
\begin{example}[A 4-dimensional example] 
Let \[ A:= \left( \begin {array}{cccc} 1&\frac{1}{x}& \frac{1}{x-1}&0
\\ \noalign{\medskip}0&1&0& \frac{1}{x-1}
\\ \noalign{\medskip}0&0&1&-\frac{1}{x}\\ \noalign{\medskip}0&0&0&1
\end {array} \right). 
\]
Note that this system is upper triangular, contrary to \eqref{redsys}. This will illustrate that our method can be equivalently applied to upper and lower triangular systems.
We obtain a Wei-Norman decomposition $A=M_1 + \frac1x M_2 + \frac{1}{x-1} M_3$, where
\[ M_1=\left( \begin {array}{cccc} 
1&0&0&0\\ \noalign{\medskip}0&1&0&0\\ \noalign{\medskip}0&0&1&0
\\ \noalign{\medskip}0&0&0&1\end {array} \right), \;
M_2:= \left(
\begin {array}{cccc} 0&1&0&0\\ \noalign{\medskip}0&0&0&0
\\ \noalign{\medskip}0&0&0&-1\\ \noalign{\medskip}0&0&0&0\end {array}
 \right), 
 \; M_3:=\left( \begin {array}{cccc} 0&0&1&0\\ \noalign{\medskip}0&0
&0&1\\ \noalign{\medskip}0&0&0&0\\ \noalign{\medskip}0&0&0&0
\end {array} \right).
\]
We have \[ M_4 :=[M_2,M_3] = \left( \begin {array}{cccc} 0&0&0&2\\ \noalign{\medskip}0&0&0&0
\\ \noalign{\medskip}0&0&0&0\\ \noalign{\medskip}0&0&0&0\end {array}
 \right)  . \]
All the other brackets in $\langle M_1,M_2,M_3, M_4 \rangle$ are zero and we find that 
$\textrm{Lie}(A)=\langle M_1,M_2,M_3, M_4 \rangle$. It is solvable of depth $2$: the first derived algebra (the set of all matrices in \textrm{Lie}(A) which can be written as a Lie bracket) is $\langle M_4\rangle$ and the second derived algebra is $\{0\}$. We will continue below with this example. 
\hfill$\diamond$
\end{example}

\subsection{Linear differential systems in reduced form}
We now turn to the link between $\textrm{Lie}(A)$ and differential Galois theory, based on two important results of Kolchin and Kovacic. Proofs can be found in \cite{PuSi03a}, Proposition~1.31 and Corollary 1.32; see also \cite{BlMo10a}, Theorem~5.8, and \cite{ApCoWe13a}, {\S}~5.3 after Remark~31. Let $A\in \mathrm{Mat}(n,\mathbf{k})$; let $G$ be the differential Galois group of $[A]$ and $\mathfrak{g}$ its Lie algebra, the Galois-Lie algebra of the system $[A]$.
The first result is 
	\[ \mathfrak{g} \, \subset \, \textrm{Lie}(A).\] 
So, the Lie algebra associated to $A$, an object which is very easy to  compute, provides an ``upper bound'' on $\mathfrak{g}$. When we perform a gauge transformation $P\in \mathrm{GL}(n,\mathbf{k})$ to obtain a new system $P[A]$,  
  $G$
and $\mathfrak{g}$ are invariant while $\textrm{Lie}(P[A])$ may vary. This shows that $\mathfrak{g}$ is a lower bound on all $\textrm{Lie}(P[A])$, for all gauge transformations $P\in \mathrm{GL}(n,\mathbf{k})$. The second result of Kolchin and Kovacic is that this lower bound is reached. 
{By definition,} $\textrm{Lie}(A)$ is the Lie algebra of an algebraic connected group $H$. 
Then there exists a gauge transformation\footnote{The notation $H(\overline{\mathbf{k}})$ denotes matrices whose entries are in $\overline{\mathbf{k}}$ and satisfy all the equations defining the algebraic group $H$.}
 $P\in H(\overline{\mathbf{k}})$ such that $\mathfrak{g} = \textrm{Lie}(P[A])$. Furthermore, if $G$ is connected and under the very mild additional condition that  $\mathbf{k}$ is a $\mathcal{C}^{1}$-field\footnote{ A field $\mathbf{k}$ is a $\mathcal{C}^{1}$-field when every non-constant homogeneous polynomial $P$ over $\mathbf{k}$ has a non-trivial zero provided that the number of its variables is more than its degree. For example, $\CC(x)$ is a $\mathcal{C}^{1}$-field and any algebraic extension of a $\mathcal{C}^{1}$-field is a $\mathcal{C}^{1}$-field (Tsen's theorem).} then we may choose $P\in H({\mathbf{k}})$ (no algebraic extension). 

\begin{definition} A system $\partial Y=AY$ is in \emph{reduced form} when the Lie algebra $\textrm{Lie}(A)$ associated to $A$ is equal to the Lie algebra of the differential Galois group of $\partial Y=AY$.
\end{definition}

The results of Kolchin and Kovacic show that a reduced form always exists. We provide, in the sequel, constructive methods to obtain them when the systems are 
in block-triangular form.
They will be illustrated on our $4$-dimensional example in the next section.

\begin{example}[4-dimensional example, continued] 
We will show, in the next section, that $[A]$ is in reduced form. 
This system is easily integrated step by step and we find a fundamental solution matrix
\[ U  = e^x  \left( \begin {array}{cccc} 1&\ln  \left( x \right) &\ln  \left( x-1
 \right) &2\, \textrm{dilog} \left( x \right) +\ln  \left( x-1 \right) 
\ln  \left( x \right) \\ \noalign{\medskip}0&1&0&\ln  \left( x-1
 \right) \\ \noalign{\medskip}0&0&1&-\ln  \left( x \right) 
\\ \noalign{\medskip}0&0&0&1\end {array} \right)
\]
where  $\textrm{dilog}$ is defined by $\textrm{dilog}'(x)=\frac{\ln(x)}{1-x}$. We may take 
 $\textrm{dilog}(x) = \textrm{Li}_2(1-x)$. 
When computing $U$, the terms requiring one integration ($\ln(x)$ and $\ln(x-1)$) correspond to the terms in $M_2$ and $M_3$ in the Wei-Norman decomposition of $A$. The term $\textrm{dilog}$ comes from the existence of $M_4$, the Lie bracket of $M_2$ and $M_3$ in $\textrm{Lie}(A)$. 

The corresponding Galois group is a semi-direct product of  a $1$-dimensional torus (giving rise to the $\exp(x)$ in the solution) and of a vector group generated by $M_2$, $M_3$ (giving rise to the terms in $\ln$)  and $M_4$ (giving rise to the $\textrm{dilog}$). Its Lie algebra is $\textrm{Lie}(A)$. 
\hfill$\diamond$
\end{example}

The ideas behind this notion of reduced form have been used for inverse problems in differential Galois theory: given an algebraic group $G$, construct a differential  system $[A]$ having $G$ as its differential Galois group. It is also a technique known in differential geometry. 
Its use for direct problems in differential Galois theory is more recent. A remark in  \cite{PuSi03a} suggests that this would be a good idea. In the context of Lie-Vessiot systems, Blazquez and Morales exploit this idea in \cite{BlMo10a}. It is developed in \cite{ApDrWe16a,ApWe12b,ApWe11a} in order to study variational equations in the context of integrability of Hamiltonian systems and the Morales-Ramis-Sim\'o theory (and later in \cite{CaWe18a} to study algebraic properties of Painlev\'e equations). For irreducible systems (or systems in block diagonal form), a  criterion for reduced forms is established in \cite{ApCoWe13a} with a decision procedure. Another, much more efficient approach is given in \cite{BaClDiWe16a,BaClDiWe20a} together with generalizations of the criterion of  \cite{ApCoWe13a}.
\\
The approach described below allows, given the above results, to compute a reduced form of a block triangular linear differential system (the last case remaining after all the above contributions). It is based upon these works, notably \cite{CaWe18a}, and is constructed in \cite{DrWe19a}. 

\section{How to Compute a Reduced Form of a Reducible System}

Assume now that $\mathbf{k}=\CC(x)$, where $\CC$ is algebraically closed of characteristic zero where the derivation $\partial$ acts trivially.  It is in particular a $\mathcal{C}^1$-field. We consider a \emph{block triangular} system over the differential field $\CC (x)$ in the same form as \eqref{redsys}, that is 
$$
[A] : \, \partial Y=A Y, 
\, \textrm{ with } 
	A=\left(\begin{array}{c|c}  {A_{1}} & 0 \\\hline {S} &  {A_{2}}\end{array}\right) \in \mathrm{Mat}(n,\mathbf{k}).
$$
Let $\adiag:=  \left(\begin{array}{c|c}  {A_{1}} & 0 \\\hline 0 &  {A_{2}} \end{array}\right)$.
In what follows, we will assume that the block diagonal part ${A_{\mathrm{diag}} }$ is in {\emph{reduced form}} and we will show how to find a gauge transformation $P$ such that $P[A]$ is in reduced form. 
By \cite{DrWe19a}, Lemma 2.7, the differential Galois group is connected 
and the reduction matrix we are looking for has coefficients in $\mathbf{k}$. 
Instead of reproving all the theory (which, in this case, can be mostly found in \cite{DrWe19a}), we will work out in details a simple example where most of the required algorithmic elements appear. 
This may help convince the reader of how the method works (the details in \cite{DrWe19a} may be technical, at least in a first  {reading}). 

\subsection{Shape of the gauge transformation}
{ Let $\mathfrak{h}_{sub}$ be the the set of \emph{off-diagonal} constant matrices of the form $\left(\begin{array}{c|c} 0 & 0 \\ \hline  {S} & 0\end{array}\right)$ 
(same sizes as in relation (\ref{redsys})). We will extend the scalars to  $\mathfrak{h}_{sub}(\mathbf{k}):= \mathfrak{h}_{sub}Ê\otimes_{\CC} \mathbf{k}$, the off-diagonal matrices with coefficients in $\mathbf{k}$.}
Our first step is that we may find a reduction matrix in a very particular shape. 

\begin{lemma}[\cite{ApDrWe16a}, Lemma 3.4]\label{lemme:shape}
There exists a gauge transformation $P\in\Big\{\mathrm{Id}+B, B\in \mathfrak{h}_{sub}(\mathbf{k})\Big\}$ such that ${\partial Y = P[A]Y}$ is in reduced form.
\end{lemma}
The following is is based on an observation from \cite{ApWe12b,ApWe11a}.
Let $P=\mathrm{Id}+B$, $B\in \mathfrak{h}_{sub}(\mathbf{k})$. 
Suppose that, for all $Q \in \left\{\mathrm{Id}+B, B\in  \mathfrak{h}_{sub}(\mathbf{k})\right\}$, we have
$\textrm{Lie}(P[A]) \subseteq \textrm{Lie}(Q[P[A]])$; then, $P[A]$ is in reduced form. In other terms, no rational gauge transformation
can turn it into a system with a smaller associated Lie algebra. 
In this case, $\textrm{Lie}(P[A])$ will be the Lie algebra of the differential Galois group and this will give us transcendence relations and algebraic relations on the solutions; this will be seen on the main example of this section. 
\par 
More generally, as we can see in \cite{DrWe19a}, Section 5, if our method can reduce a system with two diagonal blocks then we can iterate this method to obtain a reduced form of a block-triangular system with an arbitrary number of blocks on the diagonal. 
 Let us illustrate this iteration on a system with three diagonal blocks of the form
\[ \left(\begin{array}{c|c|c}
A_{1}&0&0\\ \hline
S_{2,1} &A_{2}& 0 \\ \hline
S_{3,1} &S_{3,2} &A_{3}
 \end{array}\right)
 \]
 where the block diagonal part is in reduced form (see \cite{BaClDiWe20a,BaClDiWe16a} for this). We will first reduce the south-east part (which is of the same form as \eqref{redsys})
 into a form 
\[ \left(\begin{array}{c|c}
A_{2}& 0 \\ \hline
S  &A_{3}
 \end{array}\right)
 \]
Let $P_1$ be the reduction matrix. By \cite{DrWe19a}, Lemma 5.1, the following system is automatically in reduced form
 $$A_{d}:= \left(\begin{array}{c|c|c}
A_{1}&0&0\\ \hline
0 &A_{2}& 0 \\ \hline
0 &S&A_{3}
 \end{array}\right).$$
Now we perform the gauge transformation $
 \left(\begin{array}{c|c}
\mathrm{Id}& 0 \\ \hline
0  &P_1
 \end{array}\right)$
to obtain a system of the form
 \[ \left(\begin{array}{c|c|c}
A_{1}&0&0\\ \hline
\tilde{S}_{2,1} &A_{2}& 0 \\ \hline
\tilde{S}_{3,1} &S  &A_{3}
 \end{array}\right)
 \]
 (the $S_{i,j}$  may  have changed after the first reduction step). We now see that this system in the same form as \eqref{redsys} with $A_d$ as the block diagonal matrix. So a second reduction of a two-blocks triangular system allows to reduce the initial three-blocks triangular system.\par 
  This iteration is well seen in our $4$-dimensional example below.

\begin{example}[4-dimensional example, continued] 
Let \[ A:=\left( \begin {array}{cccc} 1&0&\frac{1}{x}&0\\ \noalign{\medskip}
 \frac{1}{x-1}&1&0&-\frac{1}{x}\\ \noalign{\medskip}0&0&1&0
\\ \noalign{\medskip}0&0& \frac{1}{x-1}&1\end {array}
 \right). 
\]
A simple application of a factorization algorithm shows that it is reducible. Indeed, letting 
\[P:= \left( \begin {array}{cccc} 0&0&-1&0\\ \noalign{\medskip}-1&0&0&0
\\ \noalign{\medskip}0&0&0&1\\ \noalign{\medskip}0&1&0&0\end {array}
 \right),
 \]
we have 
\[ P[A] = 
 \left( \begin {array}{cc|c|c} 1&\frac{1}{x}& \frac{1}{x-1}&0
\\ 0&1&0& \frac{1}{x-1}
\\ \hline 0&0&1&-\frac{1}{x}\\ 
\hline
0&0&0&1
\end {array} \right). 
\]
This example is, of course, particularly simple. We use it to show how to apply the iteration procedure and Lemma \ref{lemme:shape} to simplify the system or prove that it cannot be simplified further. 

Since we consider an upper triangular system, we start with the ``north-west'' corner. We let 
\[ B := \begin{pmatrix} 1 & \frac1x \\ 0 & 1 \end{pmatrix}. \]
The diagonal part is in reduced form (solutions are $e^x$ and cannot be simplified using rational functions). 
The associated Lie algebra $\textrm{Lie}(B)$ has dimension 2. Reduction would imply to have dimension $1$. 
By Lemma \ref{lemme:shape}, a reduction matrix would have the form 
\[ P:=  \begin{pmatrix} 1 &   f(x) \\ 0 & 1 \end{pmatrix}. \]
The north-east coefficient of $P[B]$ is $\frac1x - f'(x)$. The coefficient $\frac1x - f'(x)$ could never be constant (the equation $f'(x)=\frac1x$ has no rational solution, the simple pole $\frac1x$ cannot be canceled by the derivative of a rational function). For any choice of $f$, $\textrm{Lie}(P[B])$ will have dimension $2$. It follows that $[B]$ is in reduced form. 
So we iterate. 
\\
We now pick a bigger matrix $B$:
\[ B:=  \left( \begin {array}{cc|c} 1&\frac{1}{x}& \frac{1}{x-1} 
\\ 0&1&0 
\\ \hline0&0&1 
\end {array} \right) 
\, \textrm{ and } \, 
B_{\textrm{diag}} : = \left( \begin {array}{cc|c} 1&\frac{1}{x}& 0 
\\ 0&1&0 
\\ \hline0&0&1 
\end {array} \right).
\]
By \cite{DrWe19a}, Lemma 5.1, and by the above calculation, we find that $[B_{\textrm{diag}}]$ is in reduced form.
Lemma \ref{lemme:shape} thus shows that a reduction matrix would have the simple form
\[ P:=  \begin{pmatrix} 1 & 0& f(x) \\ 0 & 1 & g(x) \\  0&0&1 \end{pmatrix}. \]
Now $\textrm{Lie}(B_{\textrm{diag}})$ has dimension $2$ and $\textrm{Lie}(B)$ has dimension 3. A reduction matrix should therefore map $B$ to a matrix with an associated Lie algebra of dimension $2$. 
\\
We have $P[B]_{2,3} = - g'(x)$ 
so there should exist constants $g_1,g_2$ and a rational function $g(x)$ such that 
$ - g'(x) = g_1 + g_2 \frac1x$. A necessary condition is $g_2=0$ and $g(x) = g_0 {-} g_1 x$. 
Similarly, there should exist constants $f_1,f_2$ such that $P[B]_{1,3} = f_1 + f_2 \frac1x$. 
We now plug our condition on $g$ into this relation and find that there should be a rational function $f$ such that
\[ f'(x) ={-} f_1 - g_1 + (g_0 {-} f_2) \frac1x + \frac{1}{x-1}. \]
Now, because of the pole of order $1$ at $x=1$, this can never have a rational solution (whatever the values of the unknown constants). It follows that, for any choice of $f$, $\textrm{Lie}(P[B])$ will have dimension $3$.
So $P$ is in reduced form.
\\
Note that our main ingredient here has been to look for a rational solution of an inhomogeneous linear differential  equation whose right-hand side contains parameters. There exist algorithms to compute conditions on the parameters (from the right-hand side) so that such an equation has rational solutions,  see subsection \ref{intermezzo} below or \cite{Sin91}, and this will be the key to what follows. 
\\
We continue iterating the reduction process. Now we will have 
\[
 B:= \left( \begin {array}{cccc} 1&\frac{1}{x}& \frac{1}{x-1}&0
\\ \noalign{\medskip}0&1&0& \frac{1}{x-1}
\\ \noalign{\medskip}0&0&1&-\frac{1}{x}\\ \noalign{\medskip}0&0&0&1
\end {array} \right), 
\, 
B_{\textrm{diag}} :=  \left( \begin {array}{cccc} 1&\frac{1}{x}& \frac{1}{x-1}&0
\\ \noalign{\medskip}0&1&0&0\\ \noalign{\medskip}0&0&1&0
\\ \noalign{\medskip}0&0&0&1\end {array} \right)
\, \textrm{ and }Ê
P:=  \left( \begin {array}{cccc} 1&0&0&f_{{1}} \left( x \right) 
\\ \noalign{\medskip}0&1&0&f_{{2}} \left( x \right) 
\\ \noalign{\medskip}0&0&1&f_{{3}} \left( x \right) 
\\ \noalign{\medskip}0&0&0&1\end {array} \right).
\]
Using again \cite{DrWe19a}, Lemma 5.1, we see that $B_{\textrm{diag}}$ is in reduced form. 
Furthermore, $\textrm{Lie}(B_{\textrm{diag}})$ has dimension $3$ and $\textrm{Lie}(B)$ has dimension $4$. 
We compute $P[B]$.
\[ 
P[B] =  \left( \begin {array}{cccc} 1&\frac{1}{x}& \frac{1}{x-1}&{
\frac {f_{{2}} \left( x \right) }{x}}+{\frac {f_{{3}} \left( x \right) }{x-1}}-{\frac {\rm d}{{\rm d}x}}f_{{1}} \left( x \right) 
\\ \noalign{\medskip}0&1&0& \frac{1}{x-1}-{\frac {\rm d}{{\rm d}x}}f_{{2}} \left( x \right) \\ \noalign{\medskip}0&0&1&-\frac{1}{x}-{\frac {\rm d}{{\rm d}x}}f_{{3}} \left( x \right) 
\\ \noalign{\medskip}0&0&0&1\end {array} \right).
\]
The relation $P[B]_{3,4} = c_{1,3} + c_{2,3}\frac{1}{x} + c_{3,3} \frac{1}{x-1} $ 
gives conditions $c_{3,3}=0, c_{2,3}=-1$ and $f_{{3}} \left( x \right) =-c_{{1,3}}x+c_{{0,3}}$.
The same study on  $P[B]_{2,4}$ gives us $f_{{2}} \left( x \right) =-c_{{1,2}}x+c_{{0,2}}$. 
Without finishing with the last coefficient, we see that $\textrm{Lie}(P[B])$ contains matrices of the following forms
(respectively because of terms in $\frac1x$ and $\frac{1}{x-1}$):
\[ M_2:= \left( \begin {array}{cccc} 0&1&0&\star \\ \noalign{\medskip}0&0&0
&0\\ \noalign{\medskip}0&0&0&-1\\ \noalign{\medskip}0&0&0&0
\end {array} \right)  
\, \textrm{ and } \,
M_3 :=  \left( \begin {array}{cccc} 0&0&1& \star
\\ \noalign{\medskip}0&0&0&1\\ \noalign{\medskip}0&0&0&0
\\ \noalign{\medskip}0&0&0&0\end {array} \right)
\]
whose Lie bracket is 
\[ M_4:=[M_2, M_3] =  \left( \begin {array}{cccc} 0&0&0&2\\ \noalign{\medskip}0&0&0&0
\\ \noalign{\medskip}0&0&0&0\\ \noalign{\medskip}0&0&0&0\end {array}
 \right).
 \] 
 So we see that, whatever our future choices may be, $\textrm{Lie}(P[B])$ will contain $M_4$ and hence have dimension $4$. 
 This shows that our system cannot be reduced so it is in reduced form. 
 Furthermore, this suggests that our reduction conditions might have been stronger : 
 requiring $P[B]_{3,4} = c_{1,3} + c_{2,3}\frac{1}{x} + c_{3,3} \frac{1}{x} $ was not enough. 
 We could have imposed $P[B]_{3,4} = 0$ and $P[B]_{2,4} = 0$. Both these relations are easily seen to be impossible to fulfill with rational functions so our system is again seen to be in reduced form. 
  \hfill $\diamond$
  \end{example}
 To summarize what this example suggests: we need to ``cancel'' terms in the purely triangular part; this reduces to finding rational solutions of linear differential equations with parametrized right-hand sides. And the order of the computations matters: here, one needs to study the relations on $f_2$ and $f_3$ before studying relations on $f_1$. We will show, in the sequel, how to systematize these ideas, using an isotypical decomposition and an adapted flag structure, and how to make them algorithmic so that a computer algebra system may perform the calculations. 
 
In this example, we had seen that a fundamental solution matrix could be written using $e^x, \ln(x), \ln(x-1)$ and $\textrm{dilog(x)}$. 
As $[B]$ is in reduced form, $\textrm{Lie}(B)$ is the Lie algebra of the Galois group and it has dimension $4$. This shows that these four functions are transcendent and algebraically independent. 
So our calculation above (long but not hard) gave us a simple proof that $\textrm{dilog(x)}$ is algebraically independent of $e^x, \ln(x), \ln(x-1)$.

\subsection{The adjoint action of the diagonal}
We recall our notations so far. We have $n_i\times n_i$ matrices  $A_i$ with coefficients in $\mathbf{k}$  and 
$$A=\left(\begin{array}{c|c} A_{1} & 0 \\\hline S & A_{2}\end{array}\right){\in \mathrm{Mat}(n,\mathbf{k})},
\, 
\adiag :=  \left(\begin{array}{c|c}  {A_{1}} & 0 \\\hline 0 &  {A_{2}}\end{array}\right).$$
If we take two off-diagonal matrices $B_1$ and $B_2$ in $\hsub$, we have $B_1.B_2=0$. 
This allows two simple calculations. 
First, let $P:=\mathrm{Id} + \sum_i f_i B_i$, with ${f_i} \in \mathbf{k}$, $B_{i}\in \mathfrak{h}_{sub}$. Then
\begin{equation} P[A] = A + \sum_i {f_i} [ {\adiag,B_i}] - \sum_i \partial( {f_i})  B_i.
\end{equation}
Furthermore,  $[ {\adiag,B_i}] \in \mathfrak{h}_{sub}(\mathbf{k})$. 
These two calculations show that  reduction will be governed by the \emph{adjoint action} ${\Psi}: X\mapsto [ {\adiag,X}]$
of the block diagonal part $\adiag$ on $\mathfrak{h}_{sub}(\mathbf{k})$.
This adjoint action $\Psi$ is a linear map. Its matrix, on  the canonical basis of $\hsub$, is 
\[ \Psi=
A_2 \otimes  \mathrm{Id}_{\mathrm{n_{1}}} -  \mathrm{Id}_{\mathrm{n_{2}}}\otimes A_1^T.\]
When $\partial Y=\adiag Y$ has an abelian Lie algebra  we may easily compute a Jordan normal form of $\Psi: X\mapsto [ {\adiag,X}]$. Furthermore the eigenvalues of $\Psi$ belong to $\mathbf{k}$.
This is the idea behind \cite{ApDrWe16a}.  
In our case, we will need a more subtle structure, an isotypical decomposition into $\Psi$-invariant subspaces of $\hsub$. 


\begin{example}[An 8-dimensional example]
We consider a matrix $A$ given by
\[A  := 
\left( \begin {array}{cccc|cccc} 
	1&0&\frac{1}{x}&0&0&0&0&0
	\\  \frac{1}{x-1}&1&0&-\frac{1}{x}&0&0&0&0
	\\ 0&0&1&0&0&0&0&0
	\\ 0&0& \frac{1}{x-1} &1&0&0&0&0
	\\ \hline 
	{\star} & {\star} & {\star} & {\star} 
		&1&0&\frac{1}{x}&0 
\\{\star} & {\star} & {\star} & {\star}
		& \frac{1}{x-1}&1&0&-\frac{1}{x}
\\ {\star} & {\star} & {\star} & {\star}
		&0&0&1&0
\\ {\star} & {\star} & {\star} & {\star}
	&0&0&\frac{1}{x-1}&1\end {array} \right).
\]
The block diagonal part is given by two copies of our $4$-dimensional example $A_1$ (here, $A_2=A_1$) and we have shown that it was in reduced form. 
The off-diagonal part is given by 
 $$
\left( \begin{array}{cccc}
 \frac{-3x+4}{4x^2}& \frac{x-4}{4x^2}&-\frac{1}{2(x-1)}+\frac{2x-2}{x^2}& -\frac{1}{x} \\
\frac{1}{2(x-1)}+\frac{-2x+5}{x^2} &\frac{x-4}{4x^{2}}&\frac{2}{x-1}+\frac{4}{x^2}&\frac{1}{2(x-1)}+\frac{2x-2}{x^2}\\
\frac{-1}{4(x-1)}&0&\frac{3x-4}{4x^2}&\frac{x+4}{4x^2}\\
-\frac{1}{2(x-1)}&\frac{1}{4(x-1)}&\frac{1}{2(x-1)}+\frac{2x-7}{x^2}&\frac{-x+4}{4x^2}
 \end{array}\right).
 $$

As $A_2=A_1$, the matrix of the adjoint action of the diagonal on ${\hsub}$ is a (sparse) $16\times 16$
matrix given by $\Psi=A_1\otimes \mathrm{Id}_{4} - \mathrm{Id}_{4}\otimes A_1^T$. 
\hfill$\diamond$
\end{example} 

\subsubsection{Isotypical decomposition}
Recall that $\hsub$ is the $\CC$-vector space of off-diagonal matrices.
We now show  how  the adjoint action $\Psi$ of the diagonal will govern the  reduction strategy on $\hsub$.  
\begin{definition}
A vector space $W\subset \hsub $ will be called a \emph{$\Psi$-space} if $\Psi (W)\subset W\otimes_{\CC} \mathbf{k}$. 
\end{definition}
The importance of these $\Psi$-spaces is stated in the following lemma. 
\begin{lemma}[\cite{DrWe19a}, Lemma 2.11]
Let $A:=\left(\begin{array}{c|c} A_{1} & 0 \\ \hline S & A_{2}\end{array}\right)$ and assume that  $\partial Y=A Y$ is in reduced form. Then, $\textrm{Lie}(A)\cap \hsub$ is a $\Psi$-space.
\end{lemma}
So our reduction strategy will be to try to project onto the smallest possible $\Psi$-space using rational gauge transformations.
In \cite{DrWe19a}, we provide references to  algorithms to decompose and factor into $\Psi$-spaces. This is obtained using an isotypical decomposition (eigenring methods) and a flag structure. 

\begin{lemma}[Krull-Schmidt]
 The $\CC$-vector space $\hsub$  admits a unique {isotypical decomposition}
	\[Ê\hsub =\displaystyle\bigoplus_{i=1}^{\kappa} W_{i} \]
	where 
\begin{itemize}
\item each $W_{i}$ is a $\Psi$-space;
\item $W_i \simeq \nu_{i} V_{i}$,  {a direct sum of $\nu_i$ $\Psi$-spaces that are all isomorphic to an {indecomposable} $\Psi$-space $V_{i}$  }
which admits a \emph{flag decomposition}  
\[ÊV_{i}= V_{i}^{[\mu]} \supsetneq V_{i}^{[\mu-1]}\supsetneq \cdots \supsetneq V_{i}^{[1]} \supsetneq V_{i}^{[0]}=\{0\}\]
and $V_{i}^{[j]}/V_{i}^{[j-1]}$ is a sum of 
isomorphic irreducible $\Psi$-spaces for $1\leq j\leq \mu$; 
\item  For $i \neq j$, the $\Psi$-spaces $V_{i}\subset W_{i}$ and $V_{j}\subset W_{j}$ are not isomorphic.
\end{itemize}
\end{lemma}

Once this decomposition and flag structure are computed, we perform, at each stage,  
a projection on a minimal $\Psi$-subspace in ${V_{i}^{[j]}}$.  For some vectors $b_i\in \mathbf{k}^N$ and a matrix 
 $E_{i,j}$ with coefficients in $\mathbf{k}$ (obtained by linear algebra), 
this reduces to computing  all 
tuples $\left(F, c_1,\ldots,c_s\right)$, with $F \in \mathbf{k}^N$ and $c_i$ constants, such that
\[F'= E_{i,j}\cdot F + \sum_i c_i b_i .\]

The resulting system $P[A]$ will be ``minimal'': it will be in reduced form. The proof of this result is technical and can be found in \cite{DrWe19a}. We will illustrate the process on our main example. 
\begin{example}[8-dimensional example, continued] \label{psi5}
In this example, 
${\hsub}$ decomposes as a direct sum ${\hsub}=\hlie_1 \oplus \hlie_5 \oplus \hlie_{10}$ of three indecomposable
$\Psi$-spaces.
\\
We first study the adjoint action $\Psi=[{\adiag,\bullet}]$ of $\adiag$ on $\hlie_5$.
We find (see \cite{DrWe20a}) an adapted basis given by off-diagonal matrices $ {N}_{2},\dots,{N}_{6}$ with south-west blocks
\[\scriptstyle{ \tiny{
\left[ \begin {array}{cccc} 0&0&0&0\\ \noalign{\medskip}
2&0&0&0\\ \noalign{\medskip}0&0&0&0\\ \noalign{\medskip}0&0&-2&0
\end {array} \right] , 
\;
\left[ \begin {array}{cccc} 0&0&-2&0
\\ \noalign{\medskip}0&0&0&-2\\ \noalign{\medskip}0&0&0&0
\\ \noalign{\medskip}0&0&0&0\end {array} \right] , 
\;
\left[ 
\begin {array}{cccc} 1&0&0&0\\ \noalign{\medskip}0&-1&0&0
\\ \noalign{\medskip}0&0&-1&0\\ \noalign{\medskip}0&0&0&1\end {array}
 \right] , 
 \;
 \left[ \begin {array}{cccc} 0&0&0&0\\ \noalign{\medskip}0&0
&0&0\\ \noalign{\medskip}1&0&0&0\\ \noalign{\medskip}0&1&0&0
\end {array} \right] , 
\;
\left[ \begin {array}{cccc} 0&-1&0&0
\\ \noalign{\medskip}0&0&0&0\\ \noalign{\medskip}0&0&0&1
\\ \noalign{\medskip}0&0&0&0\end {array} \right].
}}
\]
The matrix of the adjoint action $\Psi$ on this basis of $ \hlie_5$ is

\[ \Psi_5 :=  \left(  \begin {array}{cc|c|cc} 
	0&0&\frac{ {1}}{x-1}&0&0
	\\ 0&0&\frac{{1}}{x}&0&0
	\\ \hline  0&0&0& \frac{{1}}{x} &\frac{{1}}{x-1}
	\\  \hline 0&0&0&0&0
	\\ 0&0&0&0&0\end {array} \right).
\]
The flag structure on $\hlie_5$  suggests the following reduction path. Try to remove elements in $\langle {N}_{5},{N}_{6}\rangle$ if possible;  then in $\langle {N}_{4}\rangle$; then in $\langle {N}_{2},{N}_{3}\rangle$. How to do this will be made clear in the next section; the flag structure guides the order in which computations should be handled.

We turn to $\hlie_{10}$. We find (see \cite{DrWe20a}) a basis adapted to the flag structure given by off-diagonal matrices ${N}_{7},\dots,{N}_{16}$ whose south-west blocks are:
\begin{eqnarray*}
&\scriptstyle{\tiny{
\left[ \begin {array}{cccc} 0&0&0&0
\\ \noalign{\medskip}0&0&2&0\\ \noalign{\medskip}0&0&0&0
\\ \noalign{\medskip}0&0&0&0\end {array} \right] , 
\;
\left[ 
\begin {array}{cccc} 0&0&0&0\\ \noalign{\medskip}-1&0&0&0
\\ \noalign{\medskip}0&0&0&0\\ \noalign{\medskip}0&0&-1&0\end {array}
 \right] , 
 \;
 \left[ \begin {array}{cccc} 0&0&1&0\\ \noalign{\medskip}0&0
&0&-1\\ \noalign{\medskip}0&0&0&0\\ \noalign{\medskip}0&0&0&0
\end {array} \right] , 
}} &
 \\ &\scriptstyle{\tiny{
\left[ \begin {array}{cccc} 0&0&0&0
\\ \noalign{\medskip}0&0&0&0\\ \noalign{\medskip}0&0&0&0
\\ \noalign{\medskip}1&0&0&0\end {array} \right] , 
\;
\left[ 
\begin {array}{cccc} \frac{1}{2}&0&0&0\\ \noalign{\medskip}0&-\frac{1}{2}&0&0
\\ \noalign{\medskip}0&0&\frac{1}{2}&0\\ \noalign{\medskip}0&0&0&-\frac{1}{2}
\end {array} \right] , 
\;
\left[ \begin {array}{cccc} \frac{1}{2}&0&0&0
\\ \noalign{\medskip}0&\frac{1}{2}&0&0\\ \noalign{\medskip}0&0&-\frac{1}{2}&0
\\ \noalign{\medskip}0&0&0&-\frac{1}{2}\end {array} \right] , 
\;
\left[ 
\begin {array}{cccc} 0&0&0&-1\\ \noalign{\medskip}0&0&0&0
\\ \noalign{\medskip}0&0&0&0\\ \noalign{\medskip}0&0&0&0\end {array}
 \right] ,}} &
 \\ &\scriptstyle{\tiny{
 \left[ \begin {array}{cccc} 0&0&0&0\\ \noalign{\medskip}0&0
&0&0\\ \noalign{\medskip}\frac{1}{2}&0&0&0\\ \noalign{\medskip}0&-\frac{1}{2}&0&0
\end {array} \right] , 
\; \left[ \begin {array}{cccc} 0&-\frac{1}{2}&0&0
\\ \noalign{\medskip}0&0&0&0\\ \noalign{\medskip}0&0&0&-\frac{1}{2}
\\ \noalign{\medskip}0&0&0&0\end {array} \right] , 
\;
\left[ 
\begin {array}{cccc} 0&0&0&0\\ \noalign{\medskip}0&0&0&0
\\ \noalign{\medskip}0&-\frac{1}{2}&0&0\\ \noalign{\medskip}0&0&0&0
\end {array} \right]. 
}}&
\end{eqnarray*}
The matrix of the adjoint action $\Psi$ on this adapted basis
${N}_{7},\dots,{N}_{16}$ is: 
\[ 
\Psi_{10}  := 
\left( \begin {array}{c|cc|cccc|cc|c}
0&\frac{1}{x}& \frac{1}{x-1}&0&0&0&0&0&0&0
\\ \hline 0&0&0&\frac{1}{x}&- \frac{1}{x-1}&0&0&0&0&0
\\ 0&0&0&0&0&-\frac{1}{x}&
 	\frac{1}{x-1}&0&0&0
 \\ \hline 0&0&0&0&0&0&0&
 	\frac{1}{x-1}&0&0
 \\ 0&0&0&0&0&0&0&0&
 	\frac{1}{x-1}&0
 \\ 0&0&0&0&0&0&0&\frac{1}{x} &0&0
 \\ 0&0&0&0&0&0&0&0&\frac{1}{x}&0
\\ \hline 0&0&0&0&0&0&0&0&0& \frac{1}{x-1}
\\ 0&0&0&0&0&0&0&0&0&\frac{1}{x}
\\ \hline 0&0 &0&0&0&0&0&0&0&0
\end {array} \right). 
\] 
\end{example}

\subsubsection{Intermezzo : Reduction and Rational Solutions}\label{intermezzo}
Before we continue, let us make a quick excursion into our main algorithmic toolbox.
 we start with a simple case. We look for a condition on  $P:=\mathrm{Id} + \left(\begin{array}{c|c} 0 & 0 \\ \hline \beta & 0\end{array}\right)$
to have
\[  A=\left(\begin{array}{c|c} A_{1} & 0 \\\hline S & A_{2}\end{array}\right) 
\longrightarrow 
P[A]=\left(\begin{array}{c|c} A_{1} & 0 \\\hline  {0} & A_{2}\end{array}\right).
\]
A simple calculation shows that $\beta$ should be a rational solution of the matrix linear differential system $\beta' =  {A_2}\beta - \beta {A_1} {+} S$. 
If we let $\texttt{vec}$ denote the operator transforming a matrix into a vector by stacking its {rows}, we find (see \cite{DrWe19a}) that 
$\texttt{vec}(\beta)' = \Psi\cdot \texttt{vec}(\beta) - \texttt{vec}(S)$,
where $\Psi$ is again the adjoint action of the diagonal defined above.
So reduction will be governed by computing rational solutions of linear differential systems. When $\mathbf{k}=\CC(x)$, a computer algebra algorithm for this task has been given by Barkatou in \cite{Ba99a}, see \cite{BaClElWe12a} for a generalization to linear partial differential systems and a Maple implementation. 

Now, our general tool (also found in the above references) will be an apparently more complicated problem.  Given a matrix $\Psi$ and vectors $b_1, \ldots, b_s$, we will look for
tuples $\left(F, c_1,\ldots,c_s\right)$, with $F \in \mathbf{k}^N$ and $c_i$ constant, such that
$F'= \Psi\cdot F + \sum_i c_i b_i$. {Such tuples } form a \emph{computable} vector space and the algorithms in \cite{Ba99a,BaClElWe12a} provide this when $\mathbf{k}=\CC(x)$. Results and algorithms for general fields $\mathbf{k}$ can be found in \cite{Sin91}.
\\

We now pick concrete coefficients to show how to perform the reduction on our $8$-dimensional example. A Maple worksheet\footnote{The reader may also find a pdf version at\\ \url{http://www.unilim.fr/pages_perso/jacques-arthur.weil/DreyfusWeilReductionExamples.pdf}} with this example and the chosen coefficients may be found at \cite{DrWe20a}.

\subsubsection{Reduction on $\hlie_5$ (8-dimensional example)}
To remove all of $\hlie_5$, it would be enough to have a rational solution  to  the  system
\[ 
Y'=  \Psi_5. Y + b \; \textrm{with} \;
b
=  \left( \begin {array}{c} 
  \frac{3}{x^2} - \frac{1}{x} 
\\ \noalign{\medskip}  \frac{1}{x^2} - \frac{1}{x}
\\ \noalign{\medskip}   \frac{1}{x^2}- \frac{1}{2 \,x}
\\ \noalign{\medskip}0\\ \noalign{\medskip}   \frac{1}{x^2}
\end {array} \right)
\]
and $\Psi_5$ is given in Example \ref{psi5} (page \pageref{psi5}).
This gives us reduction equations
\[
\begin{array}{lll}
(W^{[3]}): &  \left\{ \begin{array}{ccl}
 f'_{{3,1}}   \left( x \right) &=&
\frac{1}{x^{2}}
                              \\
 f'_{{3,2}}   \left( x \right) &=&0
\end{array}\right.
                             \\
 (W^{[2]}): &    \left\{ \begin{array}{ccl}                      
f'_{{2,1}}  \left( x \right) &=&
\,{\frac {1 }{x-1}}  f_{{3,1}}  \left( x \right)
+\,{
\frac { 1  }{x}}f_{{3,2}}   \left( x \right)
+\frac{1}{x^2}-\frac{1}{2 \,x}
\end{array}\right.
                               \\
 (W^{[1]}): &         \left\{ \begin{array}{ccl}                          
 f'_{1,1}  \left( x \right) &=&
\,{\frac {1 }{x}}  f_{{2,1}}   \left( x \right)
+\frac{1}{x^2} - \frac{1}{x}                               \\
 f'_{1,2}  \left( x \right) &=&
\,{\frac {1  }{x-1}}f_{{2,1}}  \left( x \right)
+  \frac{3}{x^2} - \frac{1}{x} 
\end{array}\right.
\end{array}
\]
The first two equations correspond to the highest level $W^{[3]}$ of the flag.
To remove an element from $W^{[3]}$, there should be a rational solution to 
the equation $y'= c_{1}.\frac{1}{x^2} + c_{2}.0$. The $\CC$-vector space of 
{pairs $(c_1,c_2)\in \CC^{2}$ such that there exists   $f\in \mathbf{k}$ with} $f' = c_{1}.\frac{1}{x^2} + c_{2}.0$  is found to be $2$-dimensional;
for $\underline{c}=(1,0)$, we have $f_{3,1}:=-\frac{1}{x}+c_{3,1}$; 
for $\underline{c}=(0,1)$, we have $f_{3,2}:= c_{3,2}$, where the $c_{3,i}$ are arbitrary constants (their importance will soon be visible). Our gauge transformation is
$P^{[3]}= \textrm{Id}  + f_{3,1} N_6 + f_{3,2} N_5 $
and $A^{[2]}:=P^{[3]}[A]$ does not contain any terms from $W^{[3]}$.
\par
Now $W^{[2]}$ is $1$-dimensional. The equation for the reduction on $W_2^{[2]}$ is now
$$\begin{array}{lll}
y'	&=& {\frac {1 }{x-1}}  f_{{3,1}}  \left( x \right)
+\,{\frac { 1  }{x}}f_{{3,2}}   \left( x \right)
+\frac{1}{x^2} - \frac{1}{2 \,x}
\\ 
	&=& \frac12{\frac {2\,c_{{3,2}}+1}{x}}+{\frac { c_{{3,1}}-1}{x-1}}+\frac {1}{{x}^{2}}. 
\end{array}	$$
We have necessary and sufficient conditions on the parameters  $c_{3,i}$
to have a rational solution, namely $c_{{3,1}}=1$, $c_{{3,2}}=-\frac12$ 
and then a general rational solution $f_{2,1}:= \frac{-1}{x} + c_{2,1}$.
Our new gauge transformation is 
	$P^{[2]}= \mathrm{Id}   + (-\frac{1}{x} + c_{2,1}) N_4$
and $A^{[1]}:=P^{[2]}[A^{[2]}]$ does not contain any term from $W^{[2]}$ any more. 
\par
Finally, we look for all $(c_1,c_2)\in \CC^{2}$ such that $c_1 f_{1,1}+c_2 f_{2,2}$ is rational: we
look for non-zero pairs $(c_1,c_2)\in \CC^{2}$ such that there exists a rational solution $f\in \mathbf{k}$  of 
\begin{eqnarray*}
y' &=&
c_{1}\, \left( {\frac {1}{
x} \left( -\frac{1}{x}+c_{2,1} \right) }+\frac{1}{x^2}-\frac{1}{x}
 \right)+
c_{2}\, \left( {\frac {1}{x-1} \left( -\frac{1}{x} +c_{2,1}
 \right) }+  \frac{3}{x^2} -\frac{1}{x} \right) 
\\
&=&
{\frac {3\,c_{2}}{{x}^{2}}}+{\frac {c_{1}\, \left( c_{2,1}-1 \right) }{x}}+{\frac {c_{2}\, 
\left( c_{2,1}-1 \right) }{x-1}}.
 \end{eqnarray*}	
This integral is rational if and only if both residues are zero. As the solution $c_1=c_2=0$ is not admissible, we see that a necessary and sufficient condition is $c_{2,1}=1$.
 The set 
of desired {pairs} $(c_1,c_2)$ is of dimension $2$. For $\underline{c}=(1,0)$, we have 
{$f_{1,1}:= c_{1,1}$, for $\underline{c}=(0,1)$, we have $f_{1,2}:=-\frac{3}{x}+c_{1,2}$}, where the $c_{1,i}$ are  constants and can be chosen arbitrarily. So our last gauge transformation matrix will be 
$P^{[1]} = \textrm{Id}  -\frac3x N_{2}$ and the reduction matrix on $\hlie_5$ is 
\[ P_5:=P^{[3]}P^{[2]}P^{[1]} 
	= \textrm{Id} 
	 -\frac{3}{x} N_{2} + \left( 1-  \frac{1}{x} \right) N_{{4}}- \frac{1}{2}\,N_{{5}}
	 + \left( 1- \frac{1}{x} \right) N_{{6}}.
%
\] 
The resulting matrix $\tilde{A}:=P_5[A]$ contains no terms from  $\hlie_5$. 

\subsubsection{Reduction on $\hlie_{10}$ (8-dimensional example)}
The matrix $\Psi_{10}$   
is given in Example \ref{psi5} page \pageref{psi5}.
The reduction equations are now
\[\begin{array}{cl}
(W^{[5]}) : &  \left\{  f_{5,1}' \left( x \right) =0 \right.
\\
(W^{[4]}):  &  \left\{\begin{array}{ccl} 
	 f_{4,1}'(x) &=& \,{\frac {1  }{x}}f_{5,1}(x)  -\frac{1}{2x}
                              \\
	 f_{4,2}'(x) &=& \,{\frac {1  }{x-1}}f_{5,1}(x) - \frac{1}{2(x-1)}
\end{array} \right.
\\
(W^{[3]}):  &  \left\{\begin{array}{ccl}
  f_{3,1}'(x) &=&
\,{\frac {1  }{x}}f_{4,1}(x) + \frac{1}{x}
                      \\
 f_{3,2}'(x) &=&
\,{\frac {1  }{x}}f_{4,2}(x)- \frac{1}{2\,x}
                         \\
 f_{3,3}'(x) &=&
\,{\frac {1 }{x-1}}f_{4,1}(x) 
                      \\
 f_{3,4}'(x) &=&
\,{\frac {1  }{x-1}}f_{4,2}(x) - \frac{1}{2\,(x-1)}
	\end{array} \right.
\\
(W^{[2]}):  &  \left\{\begin{array}{ccl}
 f_{2,1}'(x) &=&
\,{\frac {1  }{x-1}}f_{3,1}(x)-\,{
\frac {1 }{x}} f_{3,2}(x) 
	- \frac{1}{2\,(x-1)}  
                              \\
 f_{2,2}'(x) &=&
-\,{\frac {1  }{x-1}}f_{3,3}(x) + \,{
\frac { 1 }{x}}f_{3,4}(x) + \frac{1}{x^2} - \frac{1}{2\,(x-1)}
                              \end{array} \right.
                              \\
(W^{[1]}) : &  \left\{
  f_{1,1}'(x) =
\,{\frac {1 }{x-1}} f_{2,1}(x)  + 
\,{ \frac {1 }{x}} f_{2,2}(x) + \frac{2}{x^2}+\frac{1}{ (x-1)} . \right.
\end{array}
\]
We will let the reader solve this iteratively following the method from the previous section. 
This will give the following successive reductions

\includegraphics[width=2.5cm]{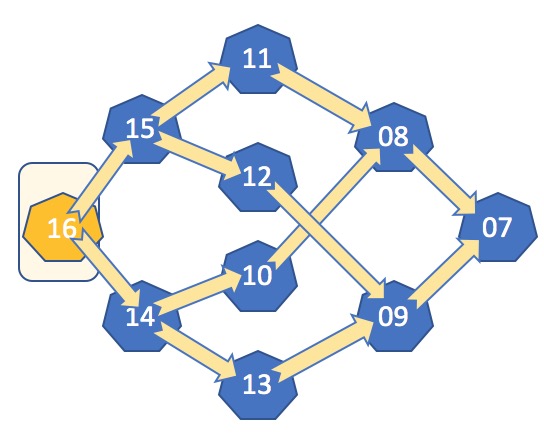} 
\includegraphics[width=2.5cm]{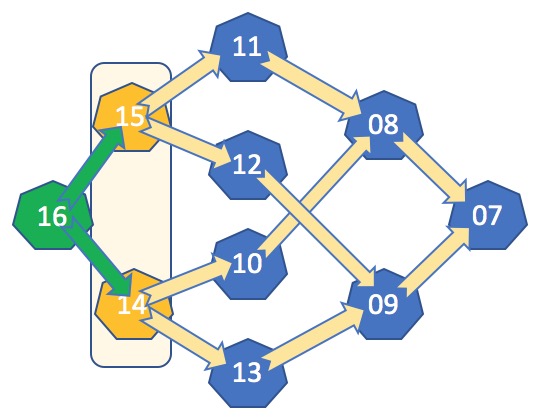} 
\includegraphics[width=2.5cm]{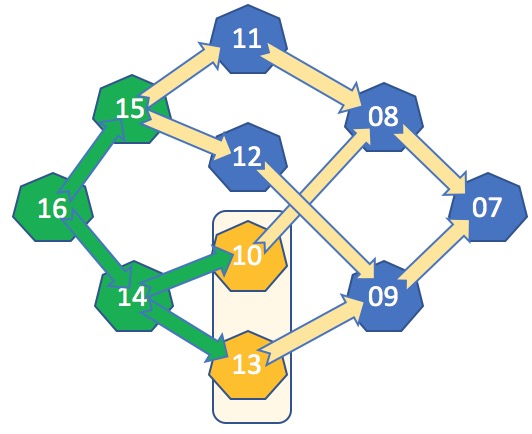} 
\includegraphics[width=2.5cm]{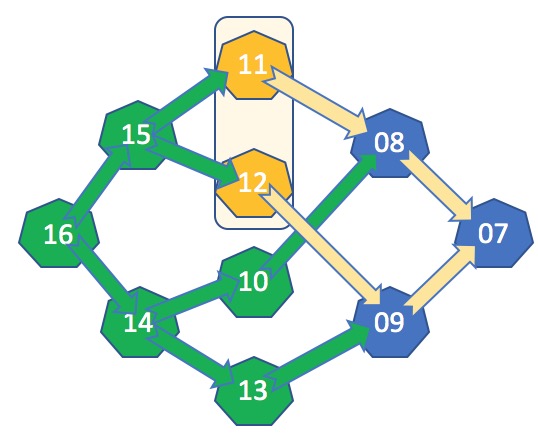} 
\\
where the green denotes parts that have been successfully removed. 
However, we reach an obstruction when trying to remove $ {N}_{11}$
(once the equation for $f_{3,1}$ has a rational solution, the  equation for $f_{3,3}(x)$ cannot have a rational solution). 	
\\
\includegraphics[width=2.5cm]{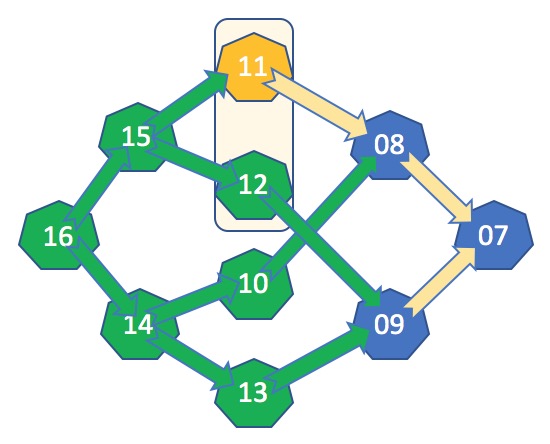} 
\includegraphics[width=2.5cm]{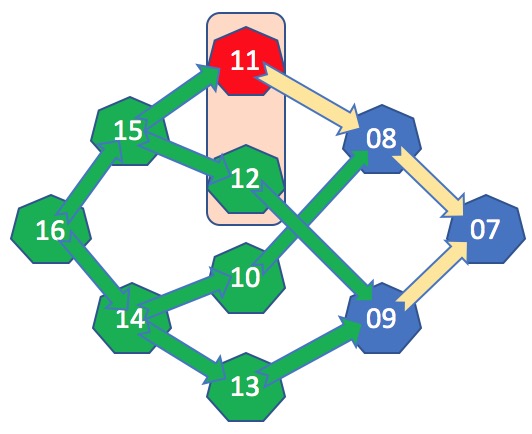}
\includegraphics[width=2.5cm]{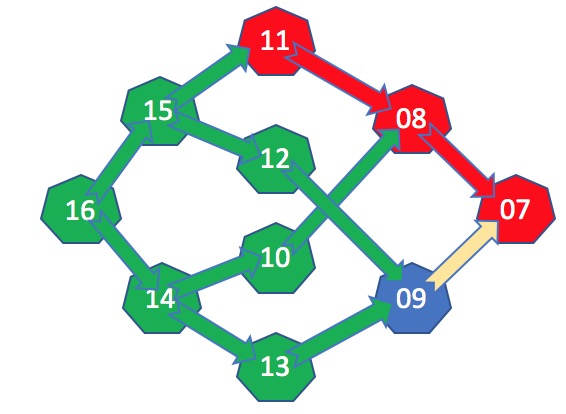}
\includegraphics[width=2.5cm]{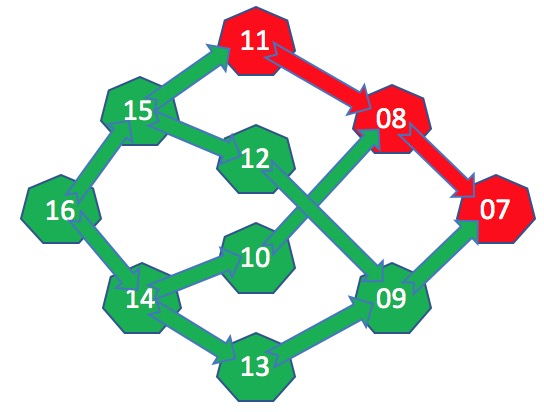} \\ 
The reduction matrix is 
\[ P_{10}:=  \textrm{Id} +
\left( c_{1,1}-\frac{1}{x} \right) N_{{7}}-{\frac {1}{x}}N_{{8}}-N_{{9}}-
\frac{1}{2}\,N_{{11}}+\frac{1}{2}\,N_{{13}}+\frac{1}{2}\,N_{{14}}-N_{{15}}+\frac{1}{2}\,N_{{16}}
\]
and we obtain 
the reduced form $A_{\textrm{red}}:=P_{10}[P_5[A]]$:
\[ 
A_{\textrm{red}} := \left( \begin {array}{cccc|cccc} 
	1&0&\frac{1}{x}&0&0&0&0&0
\\   \frac{1}{x-1}&1&0&-\frac{1}{x}&0&0&0&0
\\  0&0&1&0&0&0&0&0
\\  0&0& \frac{1}{x-1}&1&0&0&0&0
\\ \hline 
	-\frac{1}{2\,(x-1)}&0&0&0&1&0&\frac{1}{x}&0
\\  0&\frac{1}{2\,(x-1)} &0&0&          \frac{1}{x-1}&1&0&-\frac{1}{x}
\\  0&0&-\frac{1}{2\,(x-1)}&0&0&0&1&0
\\  0&0&0&\frac{1}{2\,(x-1)}&0 &0& \frac{1}{x-1}&1\end {array} \right).
\]
The associated Lie algebra is  spanned by
\begin{eqnarray*}
 {\displaystyle{\tiny \left(\begin {array}{cccc|cccc} 1&0&0&0&0&0&0&0\\  0
&1&0&0&0&0&0&0\\  0&0&1&0&0&0&0&0
\\  0&0&0&1&0&0&0&0\\\hline  0&0&0&0&1&0&0
&0\\  0&0&0&0&0&1&0&0\\  0&0&0&0&0&0
&1&0\\  0&0&0&0&0&0&0&1\end {array} \right) 
, 
\left(
\begin {array}{cccc|cccc} 0&0&0&0&0&0&0&0
\\1 &0&0&0&0&0&0&0\\  0&0&0&0&0&0&0&0
\\  0&0&1&0&0&0&0&0
\\\hline  -\frac{1}{2}&0&0&0&0&0
&0&0\\  0&\frac{1}{2}&0&0&1&0&0&0\\  0&0&-\frac{1}{2}&0&0&0&0&0
\\ 0&0&0&\frac{1}{2}&0&0&1&0\end {array}
\right) 
, 
\left(\begin {array}{cccc|cccc} 0&0
&1&0&0&0&0&0\\  0&0&0&-1&0&0&0&0\\  0
&0&0&0&0&0&0&0\\  0&0&0&0&0&0&0&0
\\\hline    0&0&0&0&0&0&1&0\\  0&0&0&0&0&0&0
&-1\\  0&0&0&0&0&0&0&0\\  0&0&0&0&0&0
&0&0\end {array} \right) }}
, 
{\displaystyle{\tiny 
\left(\begin {array}{cccc|cccc} 0&0&0&0&0&0
&0&0\\  0&0&1&0&0&0&0&0\\  0&0&0&0&0
&0&0&0\\  0&0&0&0&0&0&0&0\\\hline    0&0&{0}
&0&0&0&0&0\\  0&0&0&{0}&0&0&1&0\\  0
&0&0&0&0&0&0&0\\  0&0&0&0&0&0&0&0\end {array}
 \right) 
 , \left(\begin {array}{cccc|cccc} 0&0&0&0&0&0&0&0
\\  0&0&0&0&0&0&0&0\\  0&0&0&0&0&0&0
&0\\  0&0&0&0&0&0&0&0\\\hline    0&0&0&0&0&0
&0&0\\  0&0&\textcolor{blue}{1}&0&0&0&0&0\\  0&0&0&0&0
&0&0&0\\  0&0&0&0&0&0&0&0\end {array} \right). }}
\end{eqnarray*}
This gives us the Lie algebra $\mathfrak{g}=\textrm{Lie}( A_{\textrm{red}} )$ of the differential Galois group. 
Note that, during the reduction process, we found the two incompatible equations $f_{3,1}'(x)=
\frac {c_{4,1}+1}{x}$ and $f_{3,3}'(x)=
\frac {c_{4,1}}{x-1}$, where $c_{4,1}$ is a constant. 
There were two mutually exclusive paths: either remove $N_{11}$ or remove $N_{13}$. We removed $N_{13}$ here  {by setting $c_{4,1}=-1$}; the choice of removing $N_{11}$ ({by setting $c_{4,1}=0$}) gives a different reduced form whose associated Lie algebra is conjugated to the one we just found. {We refer to \cite{DrWe19a} for the computations in that other path.}
We also remark that two of the matrices that could not be removed from $\hsub$ are ``absorbed'' as lower triangular parts of matrices coming from $\adiag$. 
It is $5$-dimensional, whereas the Lie algebra $\textrm{Lie}(A)$ associated to the original matrix $A$ had dimension $14$.
This shows that the Picard-Vessiot extension is obtained from the Picard-Vessiot extension
$K_{\textrm{diag}}$ for $[\adiag]$ by adding only one integral and the system has indeed been transformed into a form where solving is much simpler than before - and we also have proofs of transcendence properties for the remaining objects. 

\bibliographystyle{amsalpha} 

\newcommand{\etalchar}[1]{$^{#1}$}
\providecommand{\bysame}{\leavevmode\hbox to3em{\hrulefill}\thinspace}
\providecommand{\MR}{\relax\ifhmode\unskip\space\fi MR }
\begin{thebibliography}{BCWDV16}

\bibitem[ABKM18]{AbBoKoMa18z}
Y.~Abdelaziz, S.~Boukraa, C.~Koutschan, and J.-M. Maillard, \emph{Diagonals of
  rational functions, pullbacked {$_2F_1$} hypergeometric functions and modular
  forms}, J. Phys. A \textbf{51} (2018), no.~45, 455201, 30.

\bibitem[ABKM20]{AbBoKoMa20q}
\bysame, \emph{Heun functions and diagonals of rational functions}, J. Phys. A
  \textbf{53} (2020), no.~7, 075206, 24.

\bibitem[ABRS14]{AbBlRaSc14b}
J.~Ablinger, J.~Bl\"{u}mlein, C.~G. Raab, and C.~Schneider, \emph{Iterated
  binomial sums and their associated iterated integrals}, J. Math. Phys.
  \textbf{55} (2014), no.~11, 112301, 57.

\bibitem[AMCW13]{ApCoWe13a}
A.~Aparicio-Monforte, \'E. Compoint, and J.-A. Weil, \emph{A
  characterization of reduced forms of linear differential systems}, Journal of
  Pure and Applied Algebra \textbf{217} (2013), no.~8, 1504--1516.

\bibitem[AMDW16]{ApDrWe16a}
A.~Aparicio-Monforte, T.~Dreyfus, and J.-A. Weil, \emph{Liouville
  integrability: an effective {M}orales--{R}amis--{S}im{\'o} theorem}, Journal
  of Symbolic Computation \textbf{74} (2016), 537 -- 560.

\bibitem[AMW11]{ApWe11a}
A.~Aparicio-Monforte and J.-A. Weil, \emph{A reduction method for
  higher order variational equations of {H}amiltonian systems}, Symmetries and
  Related Topics in Differential and Difference Equations, Contemporary
  Mathematics, vol. 549, Amer. Math. Soc., Providence, RI, September 2011,
  pp.~1--15.

\bibitem[AMW12]{ApWe12b}
A.~Aparicio-Monforte and J.-A. Weil, \emph{A reduced form for linear
  differential systems and its application to integrability of {H}amiltonian
  systems}, Journal of Symbolic Computation \textbf{47} (2012), no.~2, 192 --
  213.

\bibitem[AvH97]{AbHo97t}
S.~A. Abramov and M. van Hoeij, \emph{A method for the integration of
  solutions of {O}re equations}, Proceedings of the 1997 {I}nternational
  {S}ymposium on {S}ymbolic and {A}lgebraic {C}omputation ({K}ihei, {HI}), ACM,
  New York, 1997, pp.~172--175.

\bibitem[AvH99]{AbHo99b}
\bysame, \emph{Integration of solutions of linear functional equations},
  Integral Transform. Spec. Funct. \textbf{8} (1999), no.~1-2, 3--12.

\bibitem[AZ10]{AyZu10a}
M. Ayoul and Nguyen~Tien Zung, \emph{Galoisian obstructions to
  non-{H}amiltonian integrability}, C. R. Math. Acad. Sci. Paris \textbf{348}
  (2010), no.~23-24, 1323--1326.

\bibitem[Bar99]{Ba99a}
M.~A. Barkatou, \emph{On rational solutions of systems of linear
  differential equations}, J. Symbolic Comput. \textbf{28} (1999), no.~4-5,
  547--567.

\bibitem[Bax82]{Ba82a}
R.~J. Baxter, \emph{Exactly solved models in statistical mechanics},
  Academic Press Inc. [Harcourt Brace Jovanovich Publishers], London, 1982.

\bibitem[BBH{\etalchar{+}}09]{BoBoHaMaWeZe09a}
A.~Bostan, S.~Boukraa, S.~Hassani, J.-M. Maillard, J.-A. Weil, and N.~Zenine,
  \emph{Globally nilpotent differential operators and the square {I}sing
  model}, J. Phys. A \textbf{42} (2009), no.~12, 125206, 50.

\bibitem[BBH{\etalchar{+}}10]{BoBoHaMaWeZe10a}
A.~Bostan, S.~Boukraa, S.~Hassani, J.-M. Maillard, J.-A. Weil, N.~Zenine, and
  N.~Abarenkova, \emph{Renormalization, isogenies, and rational symmetries of
  differential equations}, Adv. Math. Phys. (2010), 44p.

\bibitem[BBH{\etalchar{+}}11]{BoBoHaHoMaWe11a}
A.~Bostan, S.~Boukraa, S.~Hassani, M.~van Hoeij, J.-M. Maillard, J.-A. Weil,
  and N.~Zenine, \emph{The {I}sing model: from elliptic curves to modular forms
  and {C}alabi-{Y}au equations}, J. Phys. A: Math. Theor. \textbf{44} (2011),
  no.~4, 045204, 44.

\bibitem[BCDVW20]{BaClDiWe20a}
M. Barkatou, T.~Cluzeau, L. Di~Vizio, and J.-A. Weil,
  \emph{Reduced forms of linear differential systems and the intrinsic
  {G}alois-{L}ie algebra of {K}atz}, SIGMA Symmetry Integrability Geom. Methods
  Appl. \textbf{16} (2020), Paper No. 054, 13. \MR{4112733}

\bibitem[BCEBW12]{BaClElWe12a}
M.~A. Barkatou, T.~Cluzeau, C. El~Bacha, and J.-A. Weil,
  \emph{Computing closed form solutions of integrable connections}, Proceedings
  of the 36th international symposium on Symbolic and algebraic computation
  (New York, NY, USA), ISSAC '12, ACM, 2012.

\bibitem[BCWDV16]{BaClDiWe16a}
M. Barkatou, T.~Cluzeau, J.-A. Weil, and L. Di~Vizio, \emph{Computing
  the lie algebra of the differential galois group of a linear differential
  system}, Proceedings of the ACM on International Symposium on Symbolic and
  Algebraic Computation, 2016, pp.~63--70.

\bibitem[Ber01]{Be01a}
D.~Bertrand, \emph{Unipotent radicals of differential {G}alois group and
  integrals of solutions of inhomogeneous equations}, Math. Ann. \textbf{321}
  (2001), no.~3, 645--666.

\bibitem[BHM{\etalchar{+}}07]{BoHaMaMcWeZe07a}
S.~Boukraa, S.~Hassani, J.-M. Maillard, B.~M. McCoy, J.-A. Weil, and N.~Zenine,
  \emph{Fuchs versus {P}ainlev\'e}, J. Phys. A \textbf{40} (2007), no.~42,
  12589--12605.

\bibitem[BHMW14]{BoHaMaWe14a}
S.~Boukraa, S.~Hassani, J.-M. Maillard, and J.-A. Weil, \emph{Differential
  algebra on lattice green and calabi-yau operators}, J. Phys. A: Math. Theor.
  \textbf{47} (2014), no.~9, 095203.

\bibitem[BHMW15]{BoHaMaWe15a}
S.~Boukraa, S.~Hassani, J.-M. Maillard, and J.-A. Weil, \emph{Canonical
  decomposition of irreducible linear differential operators with symplectic or
  orthogonal differential galois groups}, Journal of Physics A: Mathematical
  and Theoretical \textbf{48} (2015), no.~10, 105202.

\bibitem[BKK10]{ByKaKn10v}
V.~V. Bytev, M.~Yu. Kalmykov, and B.~A. Kniehl,
  \emph{Differential reduction of generalized hypergeometric functions from
  {F}eynman diagrams: one-variable case}, Nuclear Phys. B \textbf{836} (2010),
  no.~3, 129--170.

\bibitem[BKK13]{ByKaKn13s}
\bysame, \emph{H{YPERDIRE}, {HYPER}geometric functions {DI}fferential
  {RE}duction: {MATHEMATICA}-based packages for differential reduction of
  generalized hypergeometric functions {$_pF_{p-1},F_1,F_2,F_3,F_4$}}, Comput.
  Phys. Commun. \textbf{184} (2013), no.~10, 2332--2342.

\bibitem[BS99]{BeSi99a}
P.~H. Berman and M.~F. Singer, \emph{Calculating the {G}alois group of {$L\sb
  1(L\sb 2(y))=0$}, {$L\sb 1,L\sb 2$} completely reducible operators}, J. Pure
  Appl. Algebra \textbf{139} (1999), no.~1-3, 3--23, Effective methods in
  algebraic geometry (Saint-Malo, 1998).

\bibitem[BSMR10]{BlMo10a}
D. Bl{\'a}zquez-Sanz and J.-J. Morales-Ruiz, \emph{Differential
  {G}alois theory of algebraic {L}ie-{V}essiot systems}, Differential algebra,
  complex analysis and orthogonal polynomials, Contemp. Math., vol. 509, Amer.
  Math. Soc., Providence, RI, 2010, pp.~1--58.

\bibitem[BW03]{BoWe03a}
D. Boucher and J.-A. Weil, \emph{Application of {J}.-{J}.\
  {M}orales and {J}.-{P}.\ {R}amis' theorem to test the non-complete
  integrability of the planar three-body problem}, From combinatorics to
  dynamical systems, IRMA Lect. Math. Theor. Phys., vol.~3, de Gruyter, Berlin,
  2003, pp.~163--177.

\bibitem[CH11]{CrHa11a}
T. Crespo and Z. Hajto, \emph{Algebraic groups and differential
  {G}alois theory}, Graduate Studies in Mathematics, vol. 122, American
  Mathematical Society, Providence, RI, 2011.

\bibitem[Com12]{Co12b}
T. Combot, \emph{Non-integrability of the equal mass {$n$}-body problem
  with non-zero angular momentum}, Celestial Mech. Dynam. Astronom.
  \textbf{114} (2012), no.~4, 319--340.

\bibitem[CR91]{ChRo91a}
R.~C. Churchill and D.~L. Rod, \emph{On the determination of {Z}iglin
  monodromy groups}, SIAM J. Math. Anal. \textbf{22} (1991), no.~6, 1790--1802.

\bibitem[CW18]{CaWe18a}
G. Casale and J.-A. Weil, \emph{Galoisian methods for testing
  irreducibility of order two nonlinear differential equations}, Pacific J.
  Math. \textbf{297} (2018), no.~2, 299--337.

\bibitem[DW19]{DrWe19a}
T.~Dreyfus and J.-A. Weil, \emph{Computing the {L}ie algebra of the
  differential {G}alois group: the reducible case}, ArXiv \textbf{1904.07925}
  (2019).

\bibitem[DW20]{DrWe20a}
\bysame, \emph{Maple worksheet with the examples for this paper:
  \url{http://www.unilim.fr/pages_perso/jacques-arthur.weil/DreyfusWeilReductionExamples.mw}},
  2020.

\bibitem[FdG07]{FiGr07a}
C. Fieker and W.~A. de~Graaf, \emph{Finding integral linear dependencies
  of algebraic numbers and algebraic {L}ie algebras}, LMS J. Comput. Math.
  \textbf{10} (2007), 271--287.

\bibitem[HKMZ16]{HaKoMaZe16l}
S.~Hassani, Ch. Koutschan, J.-M. Maillard, and N.~Zenine, \emph{Lattice {G}reen
  functions: the {$d$}-dimensional face-centered cubic lattice, {$d=8, 9, 10,
  11, 12$}}, J. Phys. A \textbf{49} (2016), no.~16, 164003, 30.

\bibitem[Kal06]{Ka06g}
M.~Yu. Kalmykov, \emph{Gauss hypergeometric function: reduction,
  {$\epsilon$}-expansion for integer/half-integer parameters and {F}eynman
  diagrams}, J. High Energy Phys. (2006), no.~4, 056, 21.

\bibitem[KK11]{KaKn11a}
M.~Yu. Kalmykov and B.~A. Kniehl, \emph{Counting master integrals:
  integration by parts vs. differential reduction}, Phys. Lett. B \textbf{702}
  (2011), no.~4, 268--271.

\bibitem[KK12]{KaKn12o}
\bysame, \emph{Mellin-{B}arnes representations of {F}eynman diagrams, linear
  systems of differential equations, and polynomial solutions}, Phys. Lett. B
  \textbf{714} (2012), no.~1, 103--109.

\bibitem[KK17]{KaKn17c}
\bysame, \emph{Counting the number of master integrals for sunrise diagrams via
  the {M}ellin-{B}arnes representation}, J. High Energy Phys. (2017), no.~7,
  031, front matter+27.

\bibitem[Kou13]{Ko13a}
C. Koutschan, \emph{Lattice {G}reen functions of the higher-dimensional
  face-centered cubic lattices}, J. Phys. A \textbf{46} (2013), no.~12, 125005,
  14.

\bibitem[Kov86]{Kov86}
J.~J. Kovacic, \emph{An algorithm for solving second order linear
  homogeneous differential equations}, J. Symbolic Comput. \textbf{2} (1986),
  no.~1, 3--43.

\bibitem[LSS18]{LeSmSm18g}
R.~N. Lee, A.~V. Smirnov, and V.~A. Smirnov, \emph{Solving
  differential equations for {F}eynman integrals by expansions near singular
  points}, J. High Energy Phys. (2018), no.~3, 008, front matter+14.

\bibitem[MAB{\etalchar{+}}11]{McAsBoHaMaOr11a}
B.~M. McCoy, M.~Assis, S. Boukraa, S. Hassani, J.-M. Maillard,
  W.~P. Orrick, and N. Zenine, \emph{The saga of the {I}sing
  susceptibility}, New trends in quantum integrable systems, World Sci. Publ.,
  Hackensack, NJ, 2011, pp.~287--306.

\bibitem[Mag54]{Ma54a}
W. Magnus, \emph{On the exponential solution of differential equations for
  a linear operator}, Comm. Pure Appl. Math. \textbf{7} (1954), 649--673.

\bibitem[McC10]{Mc10a}
B.~M. McCoy, \emph{Advanced statistical mechanics}, International Series of
  Monographs on Physics, vol. 146, Oxford University Press, Oxford, 2010.

\bibitem[MM12]{McMa12a}
B.~M. McCoy and J-M. Maillard, \emph{The importance of the ising model}, Prog.
  Theor. Phys. 127 (2012), 791-817, 2012.

\bibitem[Mot18]{Mo18p}
O.~V. Motygin, \emph{On evaluation of the confluent heun functions}, 2018.

\bibitem[MR20]{Mo20s}
J.-J. Morales-Ruiz, \emph{A differential {G}alois approach to path integrals},
  J. Math. Phys. \textbf{61} (2020), no.~5, 052103, 12.

\bibitem[MRR01]{MoRa01b}
J.-J. Morales-Ruiz and J.-P. Ramis, \emph{Galoisian obstructions to
  integrability of {H}amiltonian systems. {I}, {II}}, Methods Appl. Anal.
  \textbf{8} (2001), no.~1, 33--95, 97--111.

\bibitem[MRRS07]{MRS}
J.-J. Morales-Ruiz, J.-P.  Ramis, and C. Simo, \emph{Integrability of
  {H}amiltonian systems and differential {G}alois groups of higher variational
  equations}, Ann. Sci. \'Ecole Norm. Sup. (4) \textbf{40} (2007), no.~6,
  845--884.

\bibitem[MRS09]{MoSi09a}
J.-J. Morales-Ruiz and S.~Simon, \emph{On the meromorphic non-integrability of
  some {$N$}-body problems}, Discrete Contin. Dyn. Syst. \textbf{24} (2009),
  no.~4, 1225--1273.

\bibitem[MRSS05]{MoSiSi05a}
J.-J. Morales-Ruiz, C. Sim{\'o}, and S.~Simon, \emph{Algebraic proof of the
  non-integrability of {H}ill's problem}, Ergodic Theory Dynam. Systems
  \textbf{25} (2005), no.~4, 1237--1256.

\bibitem[MS02]{MiSi02a}
C. Mitschi and M.~F. Singer, \emph{Solvable-by-finite groups as
  differential {G}alois groups}, Ann. Fac. Sci. Toulouse Math. (6) \textbf{11}
  (2002), no.~3, 403--423.

\bibitem[PPR{\etalchar{+}}10]{PPRSSW10a}
O. Pujol, J.-P. P{\'e}rez, J.-P.  Ramis, C. Sim{\'o},
  S.~Simon, and J.-A. Weil, \emph{{S}winging {A}twood {M}achine:
  experimental and numerical results, and a theoretical study}, Physica D:
  Nonlinear Phenomena \textbf{239} (2010), no.~12, 1067--1081.

\bibitem[PS03]{PuSi03a}
M. van~der Put and M.~F. Singer, \emph{Galois theory of linear differential
  equations}, Grundlehren der Mathematischen Wissenschaften [Fundamental
  Principles of Mathematical Sciences], vol. 328, Springer-Verlag, Berlin,
  2003.

\bibitem[Sal14]{Sa14a}
V. Salnikov, \emph{Effective algorithm of analysis of integrability via
  the {Z}iglin's method}, Journal of Dynamical and Control Systems \textbf{20}
  (2014), no.~4, 465--474 (English).

\bibitem[Sin91]{Sin91}
M.~F. Singer, \emph{Liouvillian solutions of linear differential equations with
  {L}iouvillian coefficients}, J. Symbolic Comput. \textbf{11} (1991), no.~3,
  251--273.

\bibitem[Sin09]{Si09a}
\bysame, \emph{Introduction to the {G}alois theory of linear differential
  equations}, Algebraic theory of differential equations, London Math. Soc.
  Lecture Note Ser., vol. 357, Cambridge Univ. Press, Cambridge, 2009,
  pp.~1--82.

\bibitem[Smi12]{Sm12e}
V.~A. Smirnov, \emph{Analytic tools for {F}eynman integrals}, Springer
  Tracts in Modern Physics, vol. 250, Springer, Heidelberg, 2012.

\bibitem[ST16]{SaTa16e}
T.~M. Sadykov and S.~Tanab\'{e}, \emph{Maximally reducible monodromy of
  bivariate hypergeometric systems}, Izv. Ross. Akad. Nauk Ser. Mat.
  \textbf{80} (2016), no.~1, 235--280.

\bibitem[Tsy01]{Ts01c}
A. Tsygvintsev, \emph{The meromorphic non-integrability of the three-body
  problem}, J. Reine Angew. Math. \textbf{537} (2001), 127--149.

\bibitem[UW96]{UlWe96a}
F.~Ulmer and J.-A. Weil, \emph{Note on {K}ovacic's algorithm}, J.
  Symbolic Comput. \textbf{22} (1996), no.~2, 179--200.

\bibitem[vHW05]{HoWe05a}
M. van Hoeij and J.-A. Weil, \emph{Solving second order differential
  equations with {K}lein's theorem}, ISSAC 2005 (Beijing), ACM, New York, 2005.

\bibitem[WN63]{WeNo63a}
J. Wei and E. Norman, \emph{Lie algebraic solution of linear
  differential equations}, J. Mathematical Phys. \textbf{4} (1963), 575--581.

\bibitem[WN64]{WeNo64a}
J.~Wei and E.~Norman, \emph{On global representations of the solutions of
  linear differential equations as a product of exponentials}, Proc. Amer.
  Math. Soc. \textbf{15} (1964), 327--334.

\bibitem[Zig82]{Zi82a}
S.~L. Ziglin, \emph{Branching of solutions and nonexistence of first integrals
  in hamiltonian mechanics. {I}}, Functional Analysis and Its Applications
  \textbf{16} (1982), no.~3, 181--189.

\end{thebibliography}
\newcommand{\etalchar}[1]{$^{#1}$}
\providecommand{\bysame}{\leavevmode\hbox to3em{\hrulefill}\thinspace}

\end{document}